\documentclass[11pt]{article}
\usepackage{geometry}
\usepackage[utf8]{inputenc}
\usepackage[T1]{fontenc}
\usepackage{amsmath}
\usepackage{amssymb}
\usepackage{mathrsfs}
\usepackage{amscd}
\usepackage{textcomp}
\usepackage{hyperref}

\usepackage{amsthm}
\theoremstyle{plain}
\newtheorem{theorem}{Theorem}[section]

\newtheorem{lemma}[theorem]{Lemma}

\newtheorem{corollary}[theorem]{Corollary}
\newtheorem{proposition}[theorem]{Proposition}

\theoremstyle{definition}

\newtheorem{definition}[theorem]{Definition}
\newtheorem{point}[theorem]{}

\newtheorem{remark}[theorem]{Remark}

\theoremstyle{remark}

\usepackage{xy}
\xyoption{all}

\newdir^{ (}{{}*!/-3pt/\dir^{(}}    
\newdir^{  }{{}*!/-3pt/\dir^{}}    
\newdir_{ (}{{}*!/-3pt/\dir_{(}}    
\newdir_{  }{{}*!/-3pt/\dir_{}}

\newcounter{zahl}

\def\theenumi{(\alph{enumi})}

\def\p@enumii{\theenumi}

\numberwithin{equation}{section}

\newcommand{\DS}{\displaystyle}

\newcommand{\SC}{\scriptstyle}
\newcommand{\SSC}{\scriptscriptstyle}
\newcommand{\sca}{\langle \,.\,,.\,\rangle}

\DeclareMathOperator{\Ad}{Ad}

\DeclareMathOperator{\End}{End}

\DeclareMathOperator{\GL}{GL}

\DeclareMathOperator{\Hom}{Hom}

\DeclareMathOperator{\Id}{Id}

\DeclareMathOperator{\Int}{Int}

\DeclareMathOperator{\Lie}{Lie}

\DeclareMathOperator{\PGL}{PGL}

\DeclareMathOperator{\SL}{SL}

\DeclareMathOperator{\SO}{SO}

\DeclareMathOperator{\Sp}{Sp}

\DeclareMathOperator{\Sym}{Sym}

\DeclareMathOperator{\Tr}{Tr}

\newcommand{\ad}{{\rm ad}}

\DeclareMathOperator{\coker}{coker}

\DeclareMathOperator{\conv}{conv}

\newcommand{\der}{{\rm der}}

\DeclareMathOperator{\diag}{diag}

\DeclareMathOperator{\id}{\,id}

\DeclareMathOperator{\ord}{ord}

\DeclareMathOperator{\tr}{tr}

\renewcommand{\phi}{\varphi}
\renewcommand{\epsilon}{\varepsilon}

\usepackage{amsfonts}
\newcommand{\BOne} {{\mathchoice{\hbox{\rm1\kern-2.7pt l\kern.9pt}}
                              {\hbox{\rm1\kern-2.7pt l\kern.9pt}}
                              {\hbox{\scriptsize\rm1\kern-2.3pt l\kern.4pt}}
                              {\hbox{\scriptsize\rm1\kern-2.4pt l\kern.5pt}}}}

\newcommand{\BC}{{\mathbb{C}}}

\newcommand{\BG}{{\mathbb{G}}}
\newcommand{\BH}{{\mathbb{H}}}

\newcommand{\BN}{{\mathbb{N}}}

\newcommand{\BQ}{{\mathbb{Q}}}
\newcommand{\BR}{{\mathbb{R}}}

\newcommand{\BZ}{{\mathbb{Z}}}

\newcommand{\CB}{{\cal{B}}}

\newcommand{\CS}{{\cal{S}}}

\newcommand{\Fa}{{\mathfrak{a}}}

\newcommand{\Fg}{{\mathfrak{g}}}
\newcommand{\Fh}{{\mathfrak{h}}}

\newcommand{\Fk}{{\mathfrak{k}}}
\newcommand{\Fl}{{\mathfrak{l}}}

\newcommand{\Fn}{{\mathfrak{n}}}
\newcommand{\Fo}{{\mathfrak{o}}}
\newcommand{\Fp}{{\mathfrak{p}}}

\newcommand{\Fs}{{\mathfrak{s}}}
\newcommand{\Ft}{{\mathfrak{t}}}

\let\setminus\smallsetminus

\newcommand{\open}{^\circ}

\newcommand{\dual}{^{\SSC\lor}}

\newcommand{\mal}{^{\SSC\times}}
\newcommand{\fdot}{{\,{\SSC\bullet}\,}}

\newcommand{\ol}[1]{{\overline{#1}}}

\usepackage{ifthen}

\newcommand{\invlim}[1][]{\ifthenelse{\equal{#1}{}}
{\DS \lim_{\longleftarrow}}
{\DS \lim_{\underset{#1}{\longleftarrow}}}
}

\newcommand{\dirlim}[1][]{\ifthenelse{\equal{#1}{}}
{\DS \lim_{\longrightarrow}}
{\DS \lim_{\underset{#1}{\longrightarrow}}}
}

\usepackage{color}
\newcounter{commentcounter}

\def\?{\ 
{\bf\color{red}???}\ 
\immediate\write16{}
\immediate\write16{Warning: There was still a question mark . . . }
\immediate\write16{}}

\long\def\forget#1{}

\long\def\addlater#1{}

\def\longto{\longrightarrow}
\def\into{\hookrightarrow}
\def\onto{\twoheadrightarrow}
\def\isoto{\stackrel{}{\mbox{\hspace{1mm}\raisebox{+1.4mm}{$\SC\sim$}\hspace{-3.5mm}$\longrightarrow$}}}
\def\longinto{\lhook\joinrel\kern-0.3em-\joinrel\kern-0.4em\rightarrow}

\newbox\mybox
\def\arrover#1{\mathrel{
       \setbox\mybox=\hbox spread 1.4em{\hfil$\scriptstyle#1$\hfil}
       \vbox{\offinterlineskip\copy\mybox
             \hbox to\wd\mybox{\rightarrowfill}}}}

\def\onto{\mbox{$\kern2pt\to\kern-8pt\to\kern2pt$}}

\begin{document}


\title{\bfseries Moment map flow on real reductive Lie groups and GIT estimates}
\author{\bfseries Christoph Böhm and Urs Hartl}
\date{}

\maketitle


\begin{abstract}
A finite dimensional real vector space carrying an action of a real reductive group possesses a moment map and a stratification as defined by Kirwan and Ness. In this article we investigate the properties of the moment map, the Kirwan-Ness stratification, and Lauret's application from a more functorial, algebraic point of view. 
\end{abstract}

\section{Introduction} \label{SectIntro}

Let $G$ be a real reductive Lie group in the sense of Borel~\cite[24.6]{Borel91}, i.e.~$G$ is (the group of $\BR$-points of) a connected reductive algebraic group over $\BR$. Let $K\subset G$ be a maximal compact subgroup and let $\Fk\subset\Fg$ be their Lie algebras. Fix a Cartan decomposition $\Fg=\Fk\oplus\Fp$. Endow $\Fg$ with an $\Ad(K)$-invariant scalar product $\sca_\Fg$, such that $\langle\Fk,\Fp\rangle_\Fg=0$ and
\[
\ad(\Fk)\;\subset\;\Fs\Fo(\Fg,\sca_\Fg)\,,\quad\text{and}\quad \ad(\Fp)\;\subset\;(\Fg,\sca_\Fg)\,.
\]
Let $\sca_\Fp$ denote the restriction of $\sca_\Fg$ to $\Fp$.

Let $\rho\colon G\to\GL(V)$ be an (algebraic) representation of $G$ on a finite dimensional $\BR$-vector space $V$ and let $\pi:=D\rho|_e\colon \Fg\to\Fg\Fl(V)$ be the induced morphism of Lie algebras. Equip $V$ with a scalar product $\sca_{_V}$ for which $\rho(g)$ is orthogonal for every $g\in K$, i.e.~$\rho(K)\subset \rm{O}(V,\sca_{_V})$, and $\pi(X)$ is symmetric, respectively skew-symmetric, if $X\in\Fp$, respectively $X\in\Fk$. In this article we investigate the $K$-equivariant moment map of the action of $G$ on $V$
\[
m\colon V\setminus\{0\}\;\longto\;\Fp
\]
implicitly  defined  by
\[
\langle m(v),X\rangle_{\Fp} \;:=\; \frac{\langle \pi(X)\cdot v,v\rangle_{_V}}{\langle v,v\rangle_{_V}}\qquad\text{for all}\quad X\,\in\,\Fp\,,
\]
The moment map
gives rise to an \emph{energy map}
\[
F\colon V\setminus\{0\}\;\longto\;\BR_{\ge0}\,,\quad v\mapsto\|m(v)\|_{\Fp}^2\,.
\]
The \emph{Kirwan-Ness stratification} yields a Morse-type disjoint union 
\[
V\setminus\{0\}\;=\;\coprod_{\beta\in\CB_{BL}}\CS_\beta
\]
into strata $\CS_\beta$ 
indexed by a finite subset $\CB_{BL}\subset\Fp$ such that
\begin{enumerate}
\item[(1)] \label{KN-Strat_1}
$\ol{\CS_\beta}\setminus\CS_\beta=\coprod_{\|\beta'\|_{\Fp}>\|\beta\|_{\Fp}}\CS_{\beta'}$, and
\item[(2)] \label{KN-Strat_2}
A vector $v\in V\setminus\{0\}$ is contained in $\CS_\beta$ if and only if the negative gradient flow of $F$ starting at $v$ converges to a critical point $v_C$ of $F$ with $m(v_C)\in K\cdot\beta$.
\end{enumerate}
This has been shown by Heinzner, Schwarz, St\"otzel~\cite{HSS08}, building on earlier work of Kirwan~\cite{Kirwan84}, Ness~\cite{Ness84}, Lauret~\cite{Lauret10}; see also \cite{Boehm-Lafuentec}.

In this article we investigate the properties of the moment map, the Kirwan-Ness stratification, and Lauret's application from a more functorial, algebraic point of view. In the analogous complex case considered by Kirwan~\cite{Kirwan84} and Ness~\cite{Ness84}, the stratification $\CS_\beta$ coincides with the algebraically defined Hesselink stratification~\cite{Hesselink79}. We show in Theorem~\ref{ThmStratum_H-BL} that this is also true in the real case; see Theorem~\ref{ThmStratum_H-BL}.

\section{Background on real reductive groups} \label{SectRealRedGps}

\begin{point}\label{PointRealRed}
There are various non-equivalent definitions of a \emph{real reductive group}; see \cite[Appendix~A]{Boehm-Lafuentec}. In this article we follow Borel~\cite[24.6]{Borel91}. We recall the definition. A \emph{linear algebraic group} over $\BR$ is a group $G$ that can be defined as a subgroup and submanifold of some $\GL_n(\BR)$ given as the zero locus of polynomial equations in the coordinates of $\GL_n(\BR)$. For a point $g\in G$ we let $T_gG$ be the \emph{tangent space} of $G$ at $g$. If $e$ is the neutral element of $G$ we call $\Fg:=\Lie G:=T_eG$ the \emph{Lie algebra} of $G$. The closed subsets in the \emph{Zariski topology} on $G$ are the subsets defined as the zero locus of polynomial equations in the coordinates of the ambient $\GL_n(\BR)$. A linear algebraic group $G$ over $R$ is \emph{reductive} if it is connected in the Zariski topology and if $\{1\}$ is the only \emph{normal unipotent algebraic} subgroup of $G$. \emph{Unipotent group} means that it can be written as a subgroup of upper triangular matrices with $1$ on the diagonal in some $\GL_n(\BR)$. Examples for real reductive groups are $\GL_n, \SL_n, \PGL_n, \SO_n, \Sp_n$ and forms like $U_n, \GL_n(\BH)$, where $\BH$ are Hamilton's quaternions, etc.

For linear algebraic groups $G$ and $G'$ over $\BR$ an (\emph{algebraic}) \emph{morphism $G\to G'$ over $\BR$} is a map given by polynomial functions in the coordinates of the ambient general linear groups with coefficients in $\BR$. It is an (\emph{algebraic}) \emph{homomorphism} if in addition it is a group homomorphism. We let $\Hom_\BR(G,G')$ be the set of algebraic homomorphisms $G\to G'$ over $\BR$. We use the notation $\BG_m$ for the multiplicative group $\GL_1(\BR)=\BR\mal=\BR\setminus\{0\}$ viewed as an algebraic group. Then $\End_\BR(\BG_m)=\BZ$ where $n\in\BZ$ corresponds to the (algebraic) endomorphism $t\mapsto t^n$. We set $X^*(G,\BR):=\Hom_\BR(G,\BG_m)$ and $X_*(G,\BR):=\Hom_\BR(\BG_m,G)$. In particular, if $A$ is a \emph{split torus}, i.e.~(algebraically) isomorphic to $\BG_m^n$ for some $n$, then $X^*(A):=\Hom_\BC(G_\BC,\BG_{m,\BC})=X^*(A,\BR)\cong\BZ^n$ and $X_*(A):=\Hom_\BC(\BG_{m,\BC},G_\BC)=X_*(A,\BR)\cong\BZ^n$, where the index $\BC$ denotes base change to $\BC$. There is a perfect pairing
\begin{equation}\label{eqn:sca_A}
\sca_A \colon X^*(A,\BR)\times X_*(A,\BR)\;\longto\;\End_\BR(\BG_m)\;=\;\BZ,\quad (\chi,\lambda)\;\longmapsto\;\chi\circ\lambda\,.
\end{equation}
\end{point}

\begin{point}\label{PointInt}
Recall that for $h\in G$ the \emph{left translation} $l_h\colon G\isoto G$, $g\mapsto h\cdot g$ and the \emph{right translation} $r_h\colon G\isoto G$, $g\mapsto g\cdot h^{-1}$ are algebraic morphisms. The inner automorphism $\Int_h\colon G\isoto G$, $g\mapsto h\cdot g\cdot h^{-1}$ is an algebraic homomorphism. It is the composition $\Int_h=l_h\circ r_h=r_h \circ l_h$. Their differentials define isomorphisms of the tangent spaces 
\begin{align*}
& D l_h|_g \colon T_gG \isoto T_{hg}G\,, \qquad D r_h|_g\colon T_gG \isoto T_{gh^{-1}}G\qquad\text{and}\\
& \Ad_h :=D\Int_h|_e=D l_h|_{h^{-1}}\circ D r_h|_e = D r_h|_h\circ D l_h|_e\colon \Fg=T_eG \isoto \Fg\,.
\end{align*}
If $G=\GL_n(\BR)\subset \BR^{n\times n}$ we can identify $T_g\GL_n(\BR)=\BR^{n\times n}$ for any $g\in G$. Then we simply have $D l_h(X) = hX$ and $D r_h(X) = X h^{-1}$ and $\Ad_h(X) = hXh^{-1}$.

If $K$ is a maximal compact (abstract) subgroup of $G$, then $K$ is a real Lie group by Cartan's theorem \cite[Part~II, \S\,V.9, Corollary to Theorem~1 on page~155]{Serre92}, but in general not an algebraic group. The right translation $l_h$ induces an isomorphism $r_h\colon K\backslash G \isoto K\backslash G$, $Kg \mapsto Kgh^{-1}$ of the homogeneous space $K\backslash G$ and an isomorphism $D r_h|_{Kg}\colon T_{Kg}(K\backslash G) \isoto T_{Kgh^{-1}}(K\backslash G)$ on tangent spaces. In particular, $D r_g|_{Kg}\colon T_{Kg}(K\backslash G)\isoto T_e(K\backslash G)=\Fk\backslash\Fg$. Composing the left translation $l_h$ on $G$ with the quotient map $G\to K\backslash G$ defines a map
\begin{equation}\label{eqn:lBarh}
G \isoto G \onto K\backslash G\,,\quad g \longmapsto hg \longmapsto Khg
\end{equation}
which only depends on the coset $\bar h=Kh\in K\backslash G$ of $h$, and which we denote $l_{\bar h}$ by abuse of notation. Beware that this map does \emph{not} factor as a map $K\backslash G \to K\backslash G$.
\end{point}

\begin{point}\label{PointRep}
An (\emph{algebraic}) \emph{representation} of a linear algebraic group $G$ over $\BR$ is an algebraic homomorphism $\rho\colon G\to\GL(V)$ over $\BR$ for a finite dimensional $\BR$-vector space $V$. We write $\pi=D\rho|_e\colon\Fg\to\Fg\Fl(V)=\End_\BR(V)$ for the differential of $\rho$ at the neutral element $e\in G$. Differentiating the morphisms $\rho\circ r_h = r_{\rho(h)}\circ \rho$ and $\rho\circ l_{h^{-1}} = l_{\rho(h)^{-1}} \circ \rho$ yields the equalities of linear maps
\begin{align}\label{eqn:DrWithrho}
D r_{\rho(h)}\circ D\rho|_h & \;=\; \pi\circ D r_h\colon\; T_h G \longto T_e\GL_n(\BR)\quad\text{and} \\[1mm]
D l_{\rho(h)^{-1}}\circ D\rho|_h & \;=\; \pi\circ D l_{h^{-1}}\colon\; T_h G \longto T_e\GL_n(\BR)\,. \nonumber
\end{align}
In particular, for every $X\in \Fg$ and $h\in G$ we have
\begin{equation} \label{eqn:piAd}
(\pi\circ \Ad_h )(X) =  (\Ad_{\rho(h)}\circ \pi)(X) = \rho(h)\circ \pi(X) \circ \rho(h^{-1}) \,.
\end{equation}
The \emph{adjoint representation} of $G$ on $\Fg$ is defined as $\Ad\colon G\to \GL(\Fg)$, $h\mapsto \Ad_h$.
\end{point}

\begin{point}\label{PointGKp}
Let $G$ be a real reductive group. Let $K\subset G$ be a maximal compact (abstract) subgroup, which hence is a real Lie group, and let $\Fk=T_eK$ be its Lie algebra. In this article we will always choose the \emph{canonical Cartan decomposition $\Fg=\Fk\oplus\Fp$ with respect to $K$} of the Lie algebra $\Fg=T_eG$ of $G$ described in \cite[24.6(a),(d)]{Borel91}. It is constructed as follows Let $Z\subset G$ be the center and let $Z_d$ be the maximal split sub-torus of $Z$. Then $\Lie Z=(\Lie Z\cap\Fk)\oplus\Lie Z_d$. Let $G^\der$ be the derived group of $G$. Then $\Fg=\Lie G^\der\oplus\Lie Z$. Recall that the \emph{Killing form} is a non-degenerate bilinear form on $\Lie G^\der$ defined by $(X,Y)\mapsto\tr(\ad(X)\circ\ad(Y))$. It defines an orthogonal complement $\Fp_0\subset\Lie G^\der$ to $(\Lie G^\der\cap \Fk)$ with $\Lie G^\der=(\Lie G^\der\cap \Fk)\oplus\Fp_0$. The Killing form is negative definite on $(\Lie G^\der\cap \Fk)$ and positive definite on $\Fp_0$. We let $\Fp:=\Fp_0\oplus\Lie Z_d$, and then $\Fg=\Fk\oplus\Fp$ is a Cartan decomposition. By \cite[24.6(a)]{Borel91} we have $G=K\cdot\exp\Fp$ and there is a \emph{Cartan involution with respect to $K$}, that is an (algebraic) group isomorphism $\Theta_K\colon G\isoto G$ with $\Theta_K^2=\id_G$ whose group of fixed points is $K$. Conversely, $\Theta_K$ induces the Cartan decomposition $\Fg=\Fk\oplus\Fp$, because $\Fp$ is the $(-1)$-eigenspace of $D\Theta_K$ and $\Fk$ is the $(+1)$-eigenspace of $D\Theta_K$. Moreover, the Lie bracket satisfies
\[
[\Fk,\Fk] \subset \Fk\,,\quad [\Fk,\Fp] \subset \Fp\,,\quad [\Fp,\Fp] \subset \Fk\,.
\]
In particular, real reductive groups $G$ in our sense satisfies the hypotheses of \cite{Boehm-Lafuentec}.
\end{point}

\begin{point}\label{PointAbelianSub}
Let $\Fa\subset\Fp$ be a maximal subalgebra which contains $\Lie Z_d$. It is automatically abelian, because $[\Fa,\Fa]\subset\Fa\cap[\Fp,\Fp]\subset\Fp\cap\Fk=(0)$. Let $A\subset G$ be the smallest algebraic subgroup whose Lie algebra contains $\Fa$. Then $A\cong\BG_m^{\dim A}=(\BR\mal)^{\dim A}$ is a maximal split torus containing $\exp(\Fa)\cong(\BR_{\SSC >0})^{\dim A}$ and $\Lie A=\Fa$ by \cite[24.6(e)]{Borel91}. Moreover, $\Fp=\bigcup_{h\in K}\Ad_h(\Fa)$ by \cite[Propositiobn~A.1]{Boehm-Lafuentec}. Next, let $M:=\{g\in G\colon \Ad(g)|_\Fa=\id_\Fa\}$, i.e.~$M$ is the centralizer $Z_G(A)$ of $A$ in $G$. Then $\Lie M=\Fa\oplus\Lie {}^{0\!} M$, where ${}^{0\!} M:=M\cap K$. If $\Fs\subset\Lie{}^{0\!} M$ is a maximal abelian subalgebra, then $\Ft:=\Fs\oplus\Fa$ is a maximal abelian subalgebra in $\Fg$, whose complexification is a Cartan subalgebra of $\Fg_\BC$; compare~\cite[\S\,2.2.4]{Wallach}. The smallest algebraic subgroup $T\subset G$ for which $\Lie T$ contains $\Ft$ is a Cartan subgroup, and hence by \cite[13.17~Corollary~2]{Borel91} a maximal torus in $G$ with maximal split sub-torus $A$ and maximal compact sub-torus $S:=\exp(\Fs)\cong(S^1)^{\dim S}$. Let $N_{K\open}(\Fa):=\{g\in K\open\colon \Ad(g)(\Fa)=\Fa\}$ and $W(\Fg,\Fa):=\{\Ad(g)|_\Fa\colon g\in N_{K\open}(\Fa)\}=N_{K\open}(\Fa)/Z_{K\open}(\Fa)$ the \emph{Weyl group} of $\Fa$. It equals the relative Weyl group $N_G(A)/Z_G(A)$ by \cite[24.6(e)]{Borel91}. The Weyl group acts simply transitively on the set of \emph{Weyl chambers}. The latter are the connected components of $\Fa':=\{X\in\Fa\colon \alpha(X)\ne0\enspace\forall\;\text{roots }\alpha\}$; compare~\cite[2.1.10~Proposition]{Wallach}.
\end{point}

For an euclidian $\BR$-vector space $(V,\sca_{_V})$ we write
\begin{equation}\label{eqn:DefSym}
\Sym(V,\sca_{_V}) \;:=\; \{\,X\in\Fg\Fl(V)\colon \langle Xv,w\rangle_{_V}=\langle v,Xw\rangle_{_V} \ \forall\,v,w\in V\,\}.
\end{equation}
In particular, if $V=\BR^n$ with the standard scalar product $\langle v,w\rangle_{_V} =v^\top w$ we write $\Sym_n(\BR):=\Sym(V,\sca_{_V}) =\{\,X\in\Fg\Fl_n\colon X^\top = X\}$.

\begin{proposition}\label{PropExistSkalProd}
Let $\rho\colon G\to\GL(V)$ be an (algebraic) representation of a real reductive group $G$ on a finite dimensional $\BR$-vector space $V$. Let $K\subset G$ be a maximal compact subgroup and let $\Fg=\Fk\oplus\Fp$ be the canonical Cartan decomposition with respect to $K$ from \S\,\ref{PointGKp} and $\Theta_K$ the corresponding Cartan involution. Then there exists a scalar product $\sca_{_V}$ on $V$ with respect to which $\rho(g)$ is orthogonal for every $g\in K$, i.e.~$\rho(K)\subset \rm{O}(V,\sca_{_V})$, and $\pi(X)$ is symmetric, respectively skew-symmetric, if $X\in\Fp$, respectively $X\in\Fk$. If $\rho$ is faithful, then $\Theta_K$ is the Cartan involution obtained from restricting the Cartan involution $g\mapsto(g^\top)^{-1}$ of $\GL(V)$ to $G$, where $g^\top$ denotes the adjoint with respect to $\sca_{_V}$. In particular, $K=G\cap \rm{O}(V,\sca_{_V})$ and $\Fp=\Fg\cap\Sym(V,\sca_{_V})$; see \eqref{eqn:DefSym}.
\end{proposition}

\begin{proof}
This was proven by Borel and Harish-Chandra~\cite[13.5~Proposition]{B+HC62}. Note that in loc.\ cit.\ the following (restrictive) definition of Cartan decomposition and Cartan involution is chosen. A Cartan involution of $\GL_n(\BR)$ is a map $g\mapsto (g^\top)^{-1}$ where $g^\top$ denotes the adjoint with respect to some scalar product on $\BR^n$. For a real reductive group $G$ a faithful algebraic representation $G\subset\GL_n(\BR)$ is chosen and a Cartan involution of $G$ is by definition the restriction to $G$ of a Cartan involution of $\GL_n(\BR)$ which maps $G$ onto itself. This definition is independent of the chosen faithful representation by \cite[\S\,13.4]{B+HC62}. The canonical Cartan involution $\Theta_K$ from \S\,\ref{PointGKp} is of this form by \cite[Proof of 13.5~Proposition]{B+HC62}.
\end{proof}

\begin{point}\label{PointScaOnLieG}
We choose a faithful representation $G\subset \GL(V)$ and a scalar product $\sca_{_V}$ on $V$ as in Proposition~\ref{PropExistSkalProd}. Via an orthonormal basis of $V$ we identify $(V,\sca_{_V})$ with $(\BR^n,\sca_{_V})$ equiped with the standard scalar product $\langle x,y\rangle_{_V}=x^\top y$. We consider the scalar product on $\Fg\Fl(V)=\Fg\Fl_n$ given by $\langle X,Y\rangle_{\Fg\Fl(V)}=\Tr(X^\top Y)$. Its restriction to $\Fg\subset\Fg\Fl_n$ will be denoted $\sca_\Fg$ and will be fixed throughout. It satisfies
\begin{equation}\label{eqn:ScaOnLieG}
\Ad(K)\;\subset\;{\rm O}(\Fg,\sca_\Fg)\,,\quad\langle \Fk,\Fp\rangle_\Fg \;=\; 0\,,\quad \ad(\Fk)\;\subset\;\Fs\Fo(\Fg,\sca_\Fg)\,,\quad \ad(\Fp)\;\subset\; \Sym(\Fg,\sca_\Fg)\,.
\end{equation}
This can be seen as follows. By Proposition~\ref{PropExistSkalProd} every element $k\in K$ is orthogonal for the standard scalar product $\sca_{_V}$ on $V=\BR^n$, i.e.~satisfies $k^\top k=1$. Therefore, $\Ad_k(X)=kXk^\top$ in $\Fg\subset\Fg\Fl_n$, and hence $\langle \Ad_k(X),\Ad_k(Y)\rangle_\Fg=\Tr(kX^\top k^\top\cdot k Yk^\top)=\Tr(X^\top Y)=\langle X,Y\rangle_\Fg$, i.e. $\Ad(K)\subset{\rm O}(\Fg,\sca_\Fg)$. This also impies $ \ad(\Fk)\subset\Fs\Fo(\Fg,\sca_\Fg)$, i.e.~$X^\top=-X$ for $X\in\Fk$. On the other hand, every element $B\in\Fp$ is symmetric with respect to the standard scalar product $\sca_{_V}$ on $\BR^n$, i.e.~satisfies $B=B^\top$. Therefore, $\ad_B(X)=BX-XB=[B,X]$ satisfies
\[
\langle \ad_B(X),Y\rangle_\Fg 
=\Tr(X^\top B^\top \cdot Y)-\Tr(B^\top X^\top \cdot Y) = \Tr(X^\top \cdot BY)-\Tr(X^\top \cdot YB) 
=\langle X,\ad_B(Y)\rangle_\Fg\,,
\]
i.e. $\ad(\Fp)\subset\Sym(\Fg,\sca_\Fg)$. Finally, if $X\in\Fk$ and $Y\in\Fp$, then $\Tr(X^\top Y)=\Tr(Y^\top X)=\Tr(X Y^\top)=\Tr(-X^\top Y)$ implies $\langle\Fk,\Fp\rangle_\Fg=0$.
\end{point}

\begin{point}\label{PointTranslateK}
Let $G$ be a real reductive group and $K\subset G$ a maximal compact subgroup. For an element $h\in G$ the conjugation action $\Int_h$ translates $K$ to $\Int_h(K)=hKh^{-1}$ and $\Fp$ to $\Ad_h(\Fp)$. If $h\in K$, then $\Int_h(K)=K$ and it follows from the definition of $\Fp$ in \S\,\ref{PointGKp} that $\Ad_h(\Fp)=\Fp$.

We sometimes denote the scalar product from Proposition~\ref{PropExistSkalProd} by $\bar g=\sca_{_V}$. The conjugation $\Int_h$ changes $\bar g=\sca_{_V}$ to the scalar product ${}^{h\!}\sca_{_V}$ on $V$ (also denoted $h\cdot \bar g$) given by
\begin{equation}\label{eqn:TranslatedScal}
(h\cdot \bar g)(x,y) \;:=\; {}^{h\!}\langle x,y\rangle_{_V} \;:=\; \langle \rho(h)^{-1}x,\rho(h)^{-1}y\rangle_{_V}
\end{equation}
for all $x,y \in V$. Namely, for $g\in K$ the element $hgh^{-1}\in hKh^{-1}$ satisfies
\[
{}^{h\!}\langle \rho(hgh^{-1}) x, \rho(hgh^{-1}) y\rangle_{_V} = \langle \rho(g)\rho(h)^{-1}x, \rho(g)\rho(h)^{-1}y \rangle_{_V} = \langle \rho(h)^{-1}x,\rho(h)^{-1}y\rangle_{_V} = {}^{h\!}\langle x,y\rangle_{_V}\,,
\]
because $\rho(g)$ is orthogonal with respect to $\sca_{_V}$. In particular, if $h\in K$ then ${}^{h\!}\sca_{_V}=\sca_{_V}$. Moreover, for $X\in\Fp$ the element $\pi(\Ad_h(X))$ satisfies
\begin{align*}
{}^{h\!}\langle \pi(\Ad_h(X)) x, y\rangle_{_V} & = \langle \rho(h)^{-1}\rho(h)\pi(X)\rho(h)^{-1}x, \rho(h)^{-1}y \rangle_{_V} = \\
& = \langle \rho(h)^{-1}x, \rho(h)^{-1}\rho(h)\pi(X)\rho(h)^{-1}y \rangle_{_V} = {}^{h\!}\langle x, \pi(\Ad_h(X)) y\rangle_{_V} \,,
\end{align*}
that is $\pi(\Ad_h(X))$ is symmetric with respect to  ${}^{h\!}\sca_{_V}$.
\end{point}

\begin{point}
If $h\colon I\to G$, $t\mapsto h(t)$ is a differentiable function on an interval $I\subset \BR$, differentiating the equation $h^{-1}(t)\cdot h(t)=1$ yields
\begin{equation}\label{eqn:DhInverse}
D l_{h^{-1}(t)}\bigl(h'(t)\bigr) + D r_{h^{-1}(t)}\bigl( \tfrac{dh^{-1}}{dt}(t)\bigr) = 0
\end{equation}
in $T_eG=\Fg$. 
\end{point}

\begin{lemma}\label{LemmaLinearODE}
Let $I\subset \BR$ be an interval with $0\in I$. Let $h_0\in G$ and let $X\colon I\to\Fg$ be a continuous function. Then the differential equations for $t\in I$
\begin{equation}\label{eqn:LinearODE1}
h'(t)= D r_{h^{-1}(t)}\bigl(X(t)\bigr) \, ,\quad (\text{respectively}\quad h'(t)= D l_{h(t)}\bigl(X(t)\bigr)\ )\,,\quad \text{with}\quad h(0) = h_0
\end{equation}
are linear, and hence have solutions for all $t\in I$.
\end{lemma}

\begin{proof}
Choosing a faithful representation $\rho\colon G\into \GL_n(\BR)$, equation~\eqref{eqn:DrWithrho} implies
\begin{align}\label{eqn:LinearODE2}
\frac{d}{dt}\rho(h(t)) & = D\rho|_{h(t)}(h'(t)) = D r_{\rho(h(t))^{-1}} \circ \pi \circ D r_{h(t)}(h'(t)) = -\pi(X(t))\cdot \rho(h(t))\,, \nonumber \\
\bigl(\ \text{or}\quad \frac{d}{dt}\rho(h(t)) & = D\rho|_{h(t)}(h'(t)) = D l_{\rho(h(t))} \circ \pi \circ D l_{h(t)^{-1}}(h'(t)) = -\rho(h(t)) \cdot \pi(X(t))\,, \quad\text{respectively }\bigr), \nonumber
\end{align}
which are linear in $\rho(h(t))$. This shows that the ODEs~\eqref{eqn:LinearODE1} are linear.
\end{proof}

\section{The moment map} \label{SectMoment}

Let $G$ be a real reductive group (\S\,\ref{PointRealRed}) and let $K\subset G$ be a maximal compact subgroup. Let $\Fg=\Fk\oplus\Fp$ be the canonical Cartan decomposition with respect to $K$ from \S\,\ref{PointGKp}. Throughout this article we fix a scalar product $\sca_\Fg$ on $\Fg$ as in \S\,\ref{PointScaOnLieG} satisfying \eqref{eqn:ScaOnLieG}. We denote by $\sca_\Fp$ the restriction of $\sca_\Fg$ to $\Fp$.

Let $\rho\colon G\to\GL(V)$ be an (algebraic) representation of $G$ on a finite dimensional $\BR$-vector space $V$ and let $\pi:=D\rho|_e\colon \Fg\to\Fg\Fl(V)$ be the induced morphism of Lie algebras. Using Proposition~\ref{PropExistSkalProd} we equip $V$ with a scalar product $\bar g=\sca_{_V}$ for which $\rho(h)$ is orthogonal for every $h\in K$, i.e.~$\rho(K)\subset \rm{O}(V,\sca_{_V})$, and $\pi(X)$ is symmetric, respectively skew-symmetric, if $X\in\Fp$, respectively $X\in\Fk$.

\begin{definition}\label{DefMomentMap}
The \emph{moment map} of the action of $G$ on $V$ with respect to the scalar product $\bar g=\sca_{_V}$ on $V$
\begin{equation}\label{eqn:DefMomentMap}
m\;:=\;m^{\bar g}\colon V\setminus\{0\}\;\longto\;\Fp
\end{equation}
is implicitly defined by
\[
\langle m(v),X\rangle_{\Fp} \;:=\; \frac{\langle \pi(X)\cdot v,v\rangle_{_V}}{\langle v,v\rangle_{_V}}\qquad\text{for all}\quad X\,\in\,\Fp\,,
\]
The moment map gives rise to an \emph{energy map}
\[
F\colon V\setminus\{0\}\;\longto\;\BR_{\ge0}\,,\quad v \;\longmapsto\;\|m(v)\|_{\Fp}^2\;:=\;\langle m(v),m(v)\rangle_\Fp\;=\;\frac{\langle \pi(m(v))\cdot v,v\rangle_{_V}}{\langle v,v\rangle_{_V}}\,.
\]
\end{definition}

Recall the metric $h\cdot\bar g = {}^{h\!}\sca_{_V}$ on $V$ from \eqref{eqn:TranslatedScal}.

\begin{definition}
For $v \in V \setminus \{0\}$ and $h \in G$
we set
\begin{align}\label{eqn:mhg}
m^{h\cdot\bar g}(v)\;:=\; \Ad_h\bigl( m(\rho(h)^{-1} v)\bigr)  \;\in\;  \Ad_h(\Fp)\,,
\end{align}
We have $m^{\bar g}(v)=m(v)$. 
\end{definition}

\begin{remark}\label{rem:MomentMapForhBarg}
(a) Notice that $m^{h\cdot\bar g}(v)$ is the moment map with respect to the metrics $h\cdot\bar g={}^{h\!}\sca_{_V}$ on $V$ and $\langle Y,X\rangle_{\Ad_h(\Fp)}:=\langle \Ad_h^{-1}(Y),\Ad_h^{-1}(X)\rangle_\Fp$ on $\Ad_h(\Fp)$, because
\[
\langle m^{h\cdot\bar g}(v),X\rangle_{\Ad_h(\Fp)} \;=\; \langle m(\rho(h)^{-1}v),\Ad_h^{-1}(X)\rangle_\Fp = \frac{\langle \pi(\Ad_h^{-1}(X))\cdot \rho(h)^{-1}v,\rho(h)^{-1}v\rangle_{_V}}{\langle \rho(h)^{-1}v,\rho(h)^{-1}v\rangle_{_V}} \;=\; \frac{{}^{h\!}\langle \pi(X)\cdot v,v\rangle_{_V}}{{}^{h\!}\langle v,v\rangle_{_V}}\,,
\]
for all $X\,\in\,\Ad_h(\Fp)$.

\medskip\noindent
(b) The data $h\cdot\bar g={}^{h\!}\sca_{_V}$ and $\Ad_h(\Fp)$ and $m^{h\cdot\bar g}(v)$ only depend on the coset $hK\in G/K$. Namely, let $h\in K$. Then ${}^{h\!}\sca_{_V}=\sca_{_V}$ and $\Ad_h(\Fp)=\Fp$ and $\sca_{\Ad_h(\Fp)}=\sca_\Fp$ by \eqref{eqn:ScaOnLieG}. Therefore, $m^{h\cdot\bar g} = m^{\bar g}$, i.e.~
\begin{equation}\label{eqn:mIsKequiv}
\Ad_h(m(v)) \;=\; m(\rho(h)v)\,.
\end{equation}
In order to work with the left homogeneous space $K\backslash G$ we invert $h$. Then the data $h^{-1}\cdot\bar g={}^{h^{-1}\!}\sca_{_V}$ and $\Ad_{h^{-1}}(\Fp)$ and $m^{h^{-1}\cdot\bar g}(v)$ only depend on the coset $Kh\in K\backslash G$. 
\end{remark}

\section{Gradient and moment map flow} \label{SectTwoFlows}

Let $G$ be a real reductive group as in \ref{PointRealRed} with Lie algebra $\Fg$. Choose a maximal compact subgroup $K\subset G$ and the canonical Cartan decomposition $\Fg=\Fk\oplus\Fp$ as in \S\,\ref{PointGKp}. Let $\rho:G\to\GL(V)$ be a (not necessarily faithful) representation of $G$ on a finite dimensional Euclidean $\BR$-vector space $(V,\sca_{_V})$ where the scalar product $\sca_{_V}$ satisfies the assertions of Proposition~\ref{PropExistSkalProd}. We sometimes write $\bar g=\sca_{_V}$ for this scalar product.

The moment map $m=m^{\bar g}\colon V\setminus\{0\}\to\Fp$ is defined as in \eqref{eqn:DefMomentMap} with respect to the fixed background metric $\bar g=\sca_{_V}$ on $V$.

\begin{definition}[Gradient flow on $V$]\label{def:odeV}
For $v_0\in V\setminus\{0\}$ we consider the dynamical system on $V \setminus \{0\}$ for $t\in\BR$
\begin{align}\label{eqn:odeV}
v'(t)= -\pi\bigl(m(v(t))\bigr) \cdot v(t) \quad \text{with}\quad v(0)=v_0\,. 
\end{align}
\end{definition}

The minus sign is just convention. Note that this is a gradient flow (up to a multiplicative constant), because
\[
\tfrac{1}{\Vert v\Vert^2}\,\langle \pi(m(v))\cdot v,v\rangle_{_V} = \langle m(v),m(v)\rangle_\Fp = \Vert m(v)\Vert^2 \,.
\]
We want to relate this gradient flow to a moment map flow on the homogeneous space $K\backslash G$ defined as follows. Recall the map $l_{\bar h}$ from \eqref{eqn:lBarh} and its differential $D l_{\bar h}$ which only depend on the coset $\bar h = Kh\in K\backslash G$ of $h$.

\begin{definition}[Moment map flow]\label{lem:odeSymRn}
Fix a vector $\bar v\in V\setminus\{0\}$. For a given element $\bar h_0=Kh_0\in K\backslash G$ we consider the dynamical system $\bar h(t)\in K\backslash G$ on the homogeneous space $K\backslash G$ for $t\in\BR$ with
\begin{align} \label{eqn:odeHomogSp}
\frac{d}{dt} \bar h(t) = - D l_{\bar h(t)}\bigl(m^{\bar h^{-1}(t)\cdot\bar g}(\bar v)\bigr) \quad \text{with}\quad \bar h(0)=\bar h_0\,.
\end{align}
It indeed only depends on the coset $\bar h(t) = Kh(t)$ by Remark~\ref{rem:MomentMapForhBarg}(b). 
\end{definition}

\begin{theorem}\label{ThmEquifOfTwoFlows}
Fix a vector $\bar v\in V\setminus\{0\}$ and an element $h_0\in G$ and let $v_0=\rho(h_0)\bar v$ and $\bar h_0=Kh_0\in K\backslash G$. Let $T\in\BR_{>0}$. Then there is a bijection between
\begin{enumerate}
\item \label{ThmEquifOfTwoFlows_1}
solutions $(v(t))_{[0,T)}$ in $V$ of the gradient flow \eqref{eqn:odeV} on $V\setminus\{0\}$
\item \label{ThmEquifOfTwoFlows_2}
and solutions $(\bar h(t))_{[0,T)}$ in $K\backslash G$ of the moment map flow \eqref{eqn:odeHomogSp}.
\end{enumerate}
\end{theorem}

\begin{proof}[Proof of Theorem~\ref{ThmEquifOfTwoFlows}.]
\ref{ThmEquifOfTwoFlows_1} $\Longrightarrow$ \ref{ThmEquifOfTwoFlows_2} \ Let $(v(t))_{t\in [0,T)}$ be a solution to \eqref{eqn:odeV} for some $T>0$. Let $h(t) \in G$ be defined by the linear differential equation
\begin{align}\label{eqn:odeh1}
h'(t)= - D r_{h^{-1}(t)}\bigl(m(v(t))\bigr) \quad \text{with}\quad h(0) = h_0 \,.
\end{align}
By Lemma~\ref{LemmaLinearODE} the solution $h(t)$ of \eqref{eqn:odeh1} exists for all $t\in [0,T)$. We claim that
\begin{equation}\label{eqn:vt}
v(t) \;=\; \rho(h(t))\cdot \rho(h_0^{-1}) \cdot v_0 \;=\; \rho(h(t))\cdot \bar v \quad\text{for all}\quad t \in [0,T)\,.
\end{equation}

To prove \eqref{eqn:vt} we set $\tilde v(t):= \rho(h^{-1}(t))\cdot v(t)$. Then, for $t_0 \in [0,T)$ we have by \eqref{eqn:piAd}, \eqref{eqn:DhInverse} and \eqref{eqn:DrWithrho}
\begin{align*}
\tilde v'(t_0) 
   &=
    \Big( \tfrac{d}{dt}\big\vert_{t=t_0}\,\rho(h^{-1}(t))\Big)\cdot v(t_0) + \rho(h^{-1}(t_0))\cdot v'(t_0)\\
   &=
    \Big( \tfrac{d}{dt}\big\vert_{t=t_0}\,\rho(h^{-1}(t))\rho(h(t_0))\bigr)\Big)\cdot \rho(h^{-1}(t_0)\big)\cdot v(t_0)
     -\rho(h^{-1}(t_0))\cdot 
      \pi(m(v(t_0)))\cdot v(t_0) \\
      &=
      \pi\big(D r_{h^{-1}(t_0)}\tfrac{dh^{-1}}{dt}(t_0)\big)\cdot \tilde v(t_0)
     - \Ad_{\rho(h^{-1}(t_0))}\circ \pi\big(m(v(t_0)\big) \cdot \tilde  v(t_0) =    0\\
      &=
      -\pi\big(Dl_{h^{-1}(t_0)}h'(t_0)\big)\cdot \tilde v(t_0) 
      - \pi\big(Dl_{h^{-1}(t_0)}\circ D r_{h^{-1}(t_0)}\big(m(v(t_0))\big)\big)\cdot \tilde  v(t_0) =    0\,.
\end{align*}
This shows that $\tilde v(t_0) = \tilde v(0)= \rho(h_0^{-1})\cdot v_0=\bar v$ is constant for all $t_0\in [0,T)$ and proves \eqref{eqn:vt}.

In particular, using \eqref{eqn:mhg} the differential equation~\eqref{eqn:odeh1} can be written as
\[
h'(t) \;=\; -D l_{h(t)} \circ \Ad_{h^{-1}(t)}\bigl( m(\rho(h(t))\bar v)\bigr) \;=\; -D l_{h(t)}\bigl(m^{h(t)^{-1}\bar g}(\bar v)\bigr)\,.
\]
This shows that $\bar h(t):=Kh(t)\in K\backslash G$ is a solution of the moment map flow \eqref{eqn:odeHomogSp}.

\medskip\noindent
\ref{ThmEquifOfTwoFlows_2} $\Longrightarrow$ \ref{ThmEquifOfTwoFlows_1} \ Suppose that $(\bar h(t))_{t \in [0,T)}$ is a solution to \eqref{eqn:odeHomogSp} with $\bar h(0)=\bar h_0$. By Remark~\ref{rem:MomentMapForhBarg}(b) it induces a map $[0,T)\to\Fg$, $t\mapsto m^{\bar h(t)^{-1}\cdot\bar g}(\bar v)$. We can thus consider the differential equation for $h(t)\in G$
\begin{align}\label{eqn:odeh2}
h'(t)= -D l_{h(t)}\bigl(m^{\bar h(t)^{-1}\cdot\bar g}(\bar v)\bigr) \,, \quad \text{with} \quad h(0)=h_0 \,.
\end{align}
By Lemma~\ref{LemmaLinearODE} this is a linear ODE, and hence has a solution on all of $[0,T)$. Let
$$
  v(t) :=\rho(h(t))\cdot \bar v\,.
$$
Then we show that $(v(t))_{t \in [0,T)}$ is a solution to \eqref{eqn:odeV} with $v(0)=v_0=\rho(h_0)\cdot \bar v$.

For all $t_0 \in [0,T)$ we have by \eqref{eqn:DhInverse}, \eqref{eqn:odeh2} and \eqref{eqn:mhg}
\begin{align*}
 v'(t_0) 
  &=
  \Big( \tfrac{d}{dt}\vert_{t=t_0}\rho\bigl(h(t)h^{-1}(t_0)\bigr)\Big)\cdot \rho(h(t_0))\cdot \bar v  \\
    &=
  \pi\bigl(D r_{h(t_0)}(h'(t_0)\bigr)\cdot v(t_0)  \\ 
  &=
    -\pi\bigl( \Ad_{h(t_0)}( m^{\bar h^{-1}(t_0)\cdot\bar g}(\bar v)) \bigr) \cdot v(t_0)\\
     &=
 - \pi\bigl(  m(\rho(h(t_0))\cdot \bar v)\bigr) \cdot v(t_0)  \\ 
   &=
 - \pi\bigl( m(v(t_0))\bigr)\cdot v(t_0)  \,.
\end{align*}
This shows the claim and proves the theorem
\end{proof}

\subsection{An example} \label{SectExampleGLn}

We let $G=\GL_n(\BR)$ and $K={\rm O}_n(\BR)=\{g\in G\colon g^\top =g^{-1}\}$ where $g^\top $ denotes the transposed matrix. These groups act on the right regular representation $W=\BR^n$ via multiplication on the left. We fix the $K$-invariant scalar product $\bar{s}=\sca_{_W}$ on $W$ given by $\langle x,y\rangle_{_W} = x^\top\cdot y$. We consider on $\Fg=\Fg\Fl(V)=\End_\BR(W)$ the scalar product $\tr (A \cdot B^\top)$. With respect to this scalar product, $\Fp=\Sym_n(\BR)=\{A\in \Fg\Fl_n\colon A^\top = A \}$ is the orthogonal complement of $\Fk=T_eK$; see \eqref{eqn:ScaOnLieG}. Let
\begin{equation}\label{eqn:Sym_n^+}
\Sym_n^+(\BR) \;:=\; \bigl\{\, A\in \Fg\Fl_n\colon A^\top = A \text{ is positive definite} \,\bigr\}\;\subset\;\Fp
\end{equation}
be the cone of scalar products $(x,y)\mapsto \bar{s}(Ax,y)=\langle A\cdot x,y\rangle_{_W}$ on $W$. 

\begin{lemma}\label{lem:ScalarProducts}
For $G=\GL_n(\BR)$ and $K={\rm O}_n(\BR)$ the map
\[
\phi\colon K\backslash G \;\isoto\; \Sym_n^+(\BR)\,, \quad \bar h \;=\; Kh \;\longmapsto\; h^{-1}\cdot\bar{s} \;=\;\bar{s}(h\cdot .\,,h\cdot .\,)\,, \quad\text{respectively}\quad A=h^\top\cdot h
\]
is well defined and a $C^\infty$-diffeomorphism.
\end{lemma}

\begin{proof}
By Remark~\ref{rem:MomentMapForhBarg}(b) the map is well defined. It is injective by construction and surjective, because every symmetric positive definite matrix $A\in\Sym_n^+(\BR)$ can be written as $A=h^\top h$ for some $h\in G=\GL_n(\BR)$, such that the columns of $h^{-1}$ form an orthonormal basis of $W$ with respect to the scalar product $\bar{s}(A x,y)$. At $\bar e=Ke\in K\backslash G$ the differential of $\phi\colon Kh\mapsto h^\top h$ is
\[
D \phi|_{\bar e}\colon \Fk\backslash\Fg=T_{\bar e}(K\backslash G) \;\isoto T_{\phi(\bar e)}\Sym_n^+(\BR) \;=\;\Sym_n(\BR)\,,\quad D \phi|_{\bar e}(X+\Fk) \;=\; X^\top+ X\,.
\]
It is well defined and injective, because $\Fk=\Fo_n=\{X\in \Fg=\Fg\Fl_n\colon X^\top + X=0\}$. It is surjective, because every matrix $Y\in \Sym_n(\BR)$ can be written as $X^\top+X$ for some $X=\tfrac{1}{2}Y\in\Fg\Fl_n$. Since $\phi$ is equivariant for the action of $g\in G$ on $K\backslash G$ via $r_g\colon Kh\mapsto Khg^{-1}$ and on $\Sym_n^+(\BR)$ via $A\mapsto (g^{-1})^\top Ag^{-1}$, this shows that the differential $D \phi
|_{\bar h}$ is an isomorphism for every $\bar h\in K\backslash G$.
\end{proof}

Now fix a representation $\rho\colon G\to\GL(V)$ as in Section~\ref{SectTwoFlows}. Then $V$ is a subrepresentation of $\bigoplus_iW^{\otimes m_i}\oplus (W\dual)^{\otimes n_i}$ for suitable $m_i,n_i\in\BN$ by \cite[Proposition~3.1(a)]{DeligneHodgeCycles}. In this way, the scalar product $\sca_{_W}$ induces a scalar product $\bar g=\sca_{_V}$ on $V$ which is $K$-invariant and satisfies $\langle \rho(h)x,y\rangle_{_V} = \langle x,\rho(h^\top)y\rangle_{_V}$ for every $h\in\GL_n(\BR)$, i.e.~$\rho(h^\top)=\rho(h)^\top$.

\begin{lemma}\label{LemmaPiAndSym+}
Let $A\in\Sym_n^+(\BR)$. Then the symmetric bilinear form on $V$ defined by
\begin{equation}\label{eqn:PiAndSym+}
{}^A\bar g(x,y) \;=\; {}^A\langle x,y\rangle_{_V} \;:=\; \langle \rho(A)\cdot x,y\rangle_{_V}
\end{equation}
is a scalar product on $V$. If $s(x,y)=\bar s(Ax,y)$ is the scalar product on $\BR^n$ associated with $A$, we also write ${}^s\bar g = {}^s\sca_{_V}:={}^A\bar g$.
\end{lemma}

\begin{proof}
This follows from the map $K\backslash G \to {\rm O}(V,\sca_{_V})\backslash \GL(V)$ from Proposition~\ref{PropExistSkalProd} together with the identifications $K\backslash G\cong \Sym_n^+(\BR)$ and ${\rm O}(V,\sca_{_V})\backslash \GL(V)\cong\{\text{scalar products on }V\}$ from Lemma~\ref{lem:ScalarProducts} by writing $A=h^\top h$ and $\rho(A)=\rho(h^\top)\rho(h)=\rho(h)^\top\rho(h)$. 
\end{proof}

\begin{definition}[Metric flow]\label{def:FlowInSym}
Fix a vector $\bar v\in V\setminus\{0\}$. For a given scalar product $s_0$ on $W$ we consider the dynamical system on $\Sym_n^+(\BR)$ for $t\in\BR$
\begin{align} \label{eqn:odemetric}
\frac{d}{dt} s(t) = -2 s(t)\bigl(m^{{}^{s(t)}\bar g}(\bar v)\cdot .\,,.\,\bigr) \quad \text{with}\quad s(0)=s_0\,.
\end{align}
\end{definition}

\begin{remark}\label{RemMetricFlow}
Under the isomorphism $s(t)=\phi\bigl(\bar h(t)\bigr)={}^{\bar h^{-1}(t)\!}\sca_{_W}$ from Lemma~\ref{lem:ScalarProducts} the metric flow~\eqref{eqn:odemetric} is compatible with the moment map flow~\eqref{eqn:odeHomogSp}. The latter takes the form
\begin{equation}\label{eqn:RemMetricFlow}
\bar h'(t) \;=\; - \bar h(t)\cdot m^{\bar h^{-1}(t)\cdot \bar g}(\bar v)
\end{equation}
for $G=\GL_n(\BR)$. This compatibility is seen as follows. The metric ${}^{s(t)}\bar g$ on $V$ equals $\bar h(t)^{-1}\cdot\bar g$ by Lemma~\ref{LemmaPiAndSym+}. The differential $\bar h'(t)$ is mapped under $D \phi_{h(t)}$ to $s'(t)$. The right hand side of \eqref{eqn:RemMetricFlow} is mapped to
\[
A \;=\; -m^{\bar h^{-1}(t)\cdot \bar g}(\bar v)^\top \cdot \bar h(t)^\top\cdot \bar h(t) - \bar h(t)^\top\cdot \bar h(t)\cdot m^{\bar h^{-1}(t)\cdot \bar g}(\bar v)\,,
\]
which corresponds to the scalar product on $W$
\begin{align*}
\langle Ax,y\rangle_{_W} & \;=\; -\bigl\langle\bar h(t)\cdot x\,,\, \bar h(t)\cdot m^{\bar h^{-1}(t)\cdot \bar g}(\bar v)\cdot y \bigr\rangle_{_W} -\bigl\langle\bar h(t)\cdot m^{\bar h^{-1}(t)\cdot \bar g}(\bar v)\cdot x\,,\, \bar h(t)\cdot y \bigr\rangle_{_W} \\
& \;=\; -2 \cdot {}^{\bar h(t)^{-1}\!}\bigl\langle m^{\bar h^{-1}(t)\cdot \bar g}(\bar v)\cdot x\,,\,y  \bigr\rangle_{_W} \\
& \;=\; -2\cdot s(t)\bigl(m^{{}^{s(t)}\bar g}(\bar v)\cdot x\,,\,y \bigr) \,.
\end{align*}
\end{remark}

\begin{corollary}
Fix a vector $\bar v\in V\setminus\{0\}$ and an element $h_0\in G=\GL_n(\BR)$ and let $v_0=\rho(h_0)\bar v$ and $\bar h_0=Kh_0\in K\backslash G$ and $s_0=h_0^{-1}\cdot\bar{s}$. Let $T\in\BR_{>0}$. Then there is a bijection between the three sets of
\begin{enumerate}
\item \label{CorEquifOfTwoFlows_1}
solutions $(v(t))_{[0,T)}$ in $V$ of the gradient flow \eqref{eqn:odeV} on $V\setminus\{0\}$,
\item \label{CorEquifOfTwoFlows_2}
solutions $(\bar h(t))_{[0,T)}$ in $K\backslash G$ of the moment map flow~\eqref{eqn:odeHomogSp},
\item \label{CorEquifOfTwoFlows_3}
and solutions $(s(t))_{[0,T)}$ in $\Sym_n^+(\BR)$ of the metric flow \eqref{eqn:odemetric}.
\end{enumerate}
\end{corollary}

\begin{proof}
This follows from Theorem~\ref{ThmEquifOfTwoFlows} and Remark~\ref{RemMetricFlow}.
\end{proof}

\section{Equality of the Kirwan-Ness and the Hesselink stratifications} \label{SectStratific} 

Also over $\BR$ the Kirwan-Ness stratification is the same as the Hesselink stratification. Over the complex numbers this was observed by Kirwan~\cite[\S\,12]{Kirwan84} and Ness~\cite[\S\,9]{Ness84}. We give more details.

Let $G$ be a real reductive group and $K\subset G$ a maximal compact subgroup. Fix the Cartan decomposition $\Fg=\Fk\oplus\Fp$ from \S\,\ref{PointGKp} and a maximal torus $T$ with maximal split subtorus $A$ as in \S\,\ref{PointAbelianSub}.

\begin{point}\label{PointHesselinksM(G)}
Hesselink defines $M(G):=M(G,\BR)$ as the set of equivalence classes of pairs $(\lambda,n)\in X_*(G,\BR)\times \BN_{\SSC >0}$ where $(\lambda,n)\sim(\lambda',n')$ if and only if $n'\cdot\lambda=n\cdot\lambda'\colon t\mapsto\lambda(t)^{n'}=\lambda'(t)^n$. Write $\tfrac{1}{n}\lambda$ for the equivalence class of $(\lambda,n)$. For the tori $T$ and $A$ we have $X_*(T,\BR)=X_*(A,\BR)\cong\BZ^{\dim A}$, because every morphism $\BG_m\to T$ factors through $A$, and hence $M(T,\BR)=M(A,\BR)=X_*(A,\BR)\otimes_\BZ \BQ$. Note that $M(A)$ is a $\BQ$-vector space, but in general $M(G)$ is not, because its lacks addition. For every $\tfrac{1}{n}\lambda\in M(G)$ the image of $\lambda$ is contained in a maximal split sub-torus $A'$. 
Since $A'$, respectively $\Lie A'$ are conjugate to $A$, respectively $\Fa$ under $G$, we obtain $M(G)=\bigcup_{g\in G}\Int_g\circ M(A)$. There is a canonical embedding
\begin{align*}
\Lie\colon\;M(G) \;\longinto\; \Fg\,, \quad \tfrac{1}{n}\lambda \;\longmapsto\; \Lie(\tfrac{1}{n}\lambda) \;:=\; D\lambda(\tfrac{1}{n})\,,
\end{align*}
which restricts on $A$ to
\begin{align}\label{eqn:M_to_p}
\Lie\colon\; M(A)=\Hom_\BR(\BG_m,A)\otimes_\BZ \BQ & \;\longinto\; \Fa\enspace\subset\enspace\Fp \\
\lambda\otimes \tfrac{1}{n} \;\,\, & \:\longmapsto\; \Lie(\tfrac{1}{n}\lambda) \;:=\; D\lambda(\tfrac{1}{n}), \nonumber
\end{align}
and induces an isomorphism $M(A)\otimes_\BQ \BR\isoto\Fa$. If $\eta=\tfrac{1}{n}\lambda\in M(G)$ and $\beta=D\lambda(\tfrac{1}{n})$ is its image then one considers the parabolic subgroup of $G$ defined by
\begin{equation}\label{eqn:ParabolicQ_beta}
Q_\eta\;:=\;\bigl\{\, g\in G\colon \exists\,\lim_{t\to 0}\lambda(t)\,g\,\lambda(t)^{-1}\text{ in }G\,\bigr\}\;=\;\bigl\{\, g\in G\colon \exists \lim_{t'\to\infty}\exp(-t'\beta)\,g\,\exp(t'\beta)\text{ in }G\,\bigr\}\;=:\;Q_\beta\,,
\end{equation}
where the equality follows from $\exp(-t'\beta)=\eta(e^{-t'})=\lambda(e^{-t'/n})$.
\end{point}

\begin{point}\label{PointNormq}
We recall the $\Ad(K)$-invariant scalar product $\sca_\Fg$ from \S\,\ref{PointScaOnLieG}. We denote by $\sca_\Fp$ its restriction to $\Fp$ and to $M(A)$ via \eqref{eqn:M_to_p}. The restriction to $M(A)$ is invariant for the action of the Weyl group $W(\Fg,\Fa)$; see \S\,\ref{PointAbelianSub}. We consider on $M(G)$ the $\Int(G)$-invariant positive definite quadratic form $q_\Fp$ defined for $\Int_g \circ \tfrac{1}{n}\lambda\in M(G)$ with $\tfrac{1}{n}\lambda\in M(A)$ and $g\in G$ by
\begin{equation}\label{eqn:QuadrFormq}
q\;:=\;q_\Fp\colon M(G) \;\longto\;\BQ\,,\quad \Int_g \circ \tfrac{1}{n} \lambda \;\longmapsto \bigl\langle D\lambda(\tfrac{1}{n}), D\lambda(\tfrac{1}{n}) \bigr\rangle_\Fp\,.
\end{equation}
If $\Fa'\subset\Fp$ is another maximal (commutative) subalgebra as in \S\,\ref{PointAbelianSub}, then $\Fa'=\Ad_h(\Fa)$ for some $h\in K$ by \cite[24.7~Proposition~1]{Borel91}. Then $\bigl\langle\Ad_h(D\lambda(\tfrac{1}{n})),\Ad_h(D\lambda(\tfrac{1}{n}))\bigr\rangle_\Fp=\langle D\lambda(\tfrac{1}{n}),D\lambda(\tfrac{1}{n})\rangle_\Fp$ by the $\Ad(K)$-invariance of $\sca_\Fp$. This shows that $q_\Fp$ only depends on $\Fp$ and not on $\Fa\subset\Fp$ and $A$. 

We claim that the scalar product $\sca_\Fp$ on $X_*(A)$ takes values in $\BZ$, and that $q$ takes values in $\BQ$. Indeed, $\sca_\Fp$ is defined via an inclusion $G\subset \GL_n(\BR)$ as the restriction of $\langle X,Y\rangle_{\Fg\Fl_n}=\Tr(X^\top Y)$. The split torus $A$ is contained in a maximal split torus of $\GL_n(\BR)$ which is conjugated into the diagonal matrices by an element $S\in\GL_n(\BR)$. If $\lambda,\lambda'\colon \BG_m\to A$, then $\Ad_S(D\lambda(1))$ and $\Ad_S(D\lambda'(1))$ are diagonal matrices with integer coefficients. Since $\Fa\subset\Fp=\Fg\cap\Sym_n(\BR)$, we have $D\lambda(1)^\top=D\lambda(1)$. Therefore,
\begin{equation*}\label{eqn:TrDlambda}
\Tr\bigl(D\lambda(1)^\top\cdot D\lambda'(1)\bigr) \;=\; \Tr\bigl(D\lambda(1)\cdot D\lambda'(1)\bigr) \;=\; \Tr\bigl(\Ad_S(D\lambda(1))\cdot \Ad_S(D\lambda'(1))\bigr) \;\in\;\BZ\,.
\end{equation*}
This shows that $q$ takes values in $\BQ$. In particular, we may take $q$ as the ``norm'' considered by Hesselink \cite[(1.3), (1.4)]{Hesselink78}.

Since $\sca_\Fp$ takes values in $\BQ$ on $M(A)$, we can define the isomorphism
\begin{equation}\label{eqn:Isomphi}
\phi\colon M(A)\;\isoto\; X^*(A)_\BQ\;:=\;X^*(A)\otimes_\BZ \BQ\,,\quad\text{by}\quad \langle \phi(\eta_1),\eta_2\rangle_{\!_A}\;=\;\langle\Lie \eta_1,\Lie \eta_2\rangle_\Fp
\end{equation}
for $\eta_i\in M(A)$. In particular $\langle\phi(\eta_1),\eta_2\rangle_{\!_A} = \langle \phi(\eta_2),\eta_1\rangle_{\!_A}$ and $q(\eta)=\langle\phi(\eta),\eta\rangle_{\!_A}$.
\end{point}

\begin{point}\label{PointMeasureInstab}
Let $V$ be an (algebraic) representation of $G$, let $x\in V$ and let $\lambda\in X_*(G,\BR)$. Hesselink~\cite[(2.2)]{Hesselink78} considers the function $f\colon \BR\mal\to V,\ t\mapsto \lambda(t)\cdot x$ and defines the \emph{measure of instability of $x$ with respect to $\lambda$} 
\[
m(x,\lambda)\;:=\; \left\{
\begin{array}{cl}
-\infty & \text{if \quad $\not\hspace{-0.25em}\exists\;\lim_{t\to 0} f(t)$,}\\[1mm]
\infty & \text{if \quad $x=0$,}\\[1mm]
0 & \text{if \quad $\exists\;\lim_{t\to 0} f(t)\ne 0$,}\\[1mm]
\ord_{t=0} f & \text{if \quad $\exists\;\lim_{t\to 0} f(t)= 0$,}
\end{array} \right.
\]
where $\ord_{t=0} f$ denotes the order of vanishing of $f$ at $t=0$. Since $m(x,n\lambda)=n\cdot m(x,\lambda)$ for every $n\in\BN_{\SSC >0}$, the definition of $m(x,\lambda)$ extends to $\tfrac{1}{n}\lambda\in M(G)$ by $m(x,\tfrac{1}{n}\lambda):=\tfrac{1}{n}m(x,\lambda)$. Then $x$ is unstable (not semi-stable) in the sense of GIT if and only if there is a $\lambda\in X_*(G,\BR)$ with $m(x,\lambda)>0$; see (the proof of) \cite[(ii) on page~53]{GIT} (together with \cite{AH-LH} which allows to work over the non-algebraically closed field $\BR$). These $x$ form the \emph{null cone} $N_G(V)$ consisting of those elements of $V$ whose orbit contains $0$ in its closure.
\end{point}

\begin{point}\label{PointOptClass}
Note that $m(x,\Int_g\circ\eta)=m(g\cdot x,\eta)$ for all $g\in G$. This yields another way to compute $m(x,\eta)$ as follows. Consider the weight space decomposition
\begin{equation} \label{eqn:WeigthSpaces}
V\;=\;\bigoplus_{\chi\in X^*(A)} V_\chi\qquad\text{with}\qquad V_\chi\;:=\;\{\,x\in V\colon a\cdot x=\chi(a)\cdot x\enspace\forall\,a\in A\,\}
\end{equation}
and the projections $p_\chi\colon V\to V_\chi$. The \emph{weights of $V$} is the set $X^*(A,V):=\{\chi\in X^*(A)\colon V_\chi\ne(0)\}$. For an $x\in V\setminus\{0\}$ Hesselink defines the \emph{state} $R(x,T):=R(x,A):=\{\chi\in X^*(A)\colon p_\chi(x)\ne0\}$. If $\eta\in M(A)$, then
\[
m(x,\eta)\;=\; \left\{
\begin{array}{cl}
-\infty & \text{if \quad $\min\{\langle \chi,\eta\rangle_{\!_A}\colon \chi\in R(x,A)\}<0$,}\\[1mm]
\min\{\langle \chi,\eta\rangle_{\!_A}\colon \chi\in R(x,A)\} & \text{else.}
\end{array} \right.
\]
Also note that if $x$ is an eigenvector for another $\tilde\eta\in X_*(A)_\BQ$, that is $\tilde\eta(t)\cdot x=t^r\cdot x$ for $r\in\BQ$, then $\langle\chi,\tilde\eta\rangle_{\!_A}=r$ for all $\chi\in R(x,A)$, and hence
\begin{equation}\label{eqn:mOfSum}
\min\{\langle \chi,\eta+\tilde\eta\rangle_{\!_A}\colon \chi\in R(x,A)\} \;=\; \min\{\langle \chi,\eta\rangle_{\!_A}\colon \chi\in R(x,A)\}+r\,.
\end{equation}
Hesselink further defines
\begin{eqnarray*}
q^*_G(x) & := & \inf\{\,q(\eta)\colon \eta\in M(G),\ m(x,\eta)\ge1\,\} \qquad \text{and}\\[2mm]
\Lambda_G(x) & := & \bigl\{\eta\in M(G)\colon m(x,\eta)\ge 1, \ q(\eta)=q^*_G(x)\,\bigr\} \quad\text{the \emph{optimal class}.}
\end{eqnarray*}
Clearly, all these definitions coincide for $x$ and $r\cdot x$ for all $r\in\BR\mal$. Intuitively, if $G=A$, then $\eta\in\Lambda_G(x)$ if and only if among all elements on the sphere $\{\eta\in M(A)\colon q(\eta)=1\}$, the normalized cocharacter $\tilde\eta:=\eta/\sqrt{q(\eta)}$ pushes $x$ to $\lim\limits_{t\to 0} \tilde\eta(t)\cdot x=0$ in the fastest way. (Except that $\tilde\eta\not\in M(A)$ if $q(\eta)$ is not a square in $\BQ$.)

A subgroup $H$, like for example $H=A$, is \emph{optimal for $x$} if $q^*_H(x)=q^*_G(x)$, i.e.~if the minimum is attained on $M(H)$. Since $M(G)=\bigcup_{g\in G}\Int_g\circ M(A)$ and $m(x,\Int_g\circ\eta)=m(g\cdot x,\eta)$ there is always an element $g\cdot x$ in the $G$-orbit of $x$ for which $A$ is optimal. If $A$ is optimal for $x$, then the optimal class can be computed as follows. Recall the isomorphism $\phi$ from \eqref{eqn:Isomphi} and define $\eta=\eta(x,A)\in M(A)$ as the unique element in the convex hull of $\phi^{-1}R(x,A)$ with minimal norm $q(\eta)$. Then $\eta\ne0$ if and only if $x$ is unstable by \cite[(3.2)~Lemma]{Hesselink78}. In this case $\Lambda_A(x)=\{\eta/q(\eta)\}$ by \cite[(4.2)~Lemma]{Hesselink78}. Moreover, \cite[Theorem~5.2]{Hesselink78} shows that there is a unique parabolic subgroup $Q\subset G$ which for every $\eta\in \Lambda_G(x)$ equals the group $Q_\eta$ from \eqref{eqn:ParabolicQ_beta}, and that $\Lambda_G(x)=\{\Int_g\circ \eta/q(\eta)\colon g\in Q\}$.
\end{point}

\begin{definition}\label{DefHessStrat}
We define \emph{Hesselink's set of stratum labels} as
\[
\CB_H(A,V)\;:=\;\bigl\{\,\eta\in M(A)\colon \exists\, \chi_1,\ldots,\chi_n\in X^*(A,V)\text{ such that }\eta=\min\nolimits_q \conv(\phi^{-1}\chi_i)\,\bigr\}\,\big/\,W(\Fg,\Fa)\,,
\]
where $\conv(\phi^{-1}\chi_i)$ is the convex hull of $\phi^{-1}(\chi_1),\ldots,\phi^{-1}(\chi_n)\in M(A)$, and $\min\nolimits_q \conv(\phi^{-1}\chi_i)$ denotes the element of $\conv(\phi^{-1}\chi_i)$ which minimizes the norm $q$.

For $\eta\in\CB_H(A,V)$ Hesselink~\cite[4.1]{Hesselink79} defines the \emph{stratum}
\[
\CS^H_\eta\;:=\; G\cdot\bigl\{\,x\in V\setminus\{0\}\colon \eta/q(\eta) \in \Lambda_G(x)\,\bigr\}.
\]
In \cite[Propositions~4.2 and 4.3]{Hesselink79} he shows that $\CS^H_\eta$ is a locally closed subset of $V$ in the Zariski topology, invariant under scaling by $\BR\mal$, and that its Zariski closure satisfies
\begin{equation}\label{eqn:HesslinkStrat}
\ol{\CS^H_\eta} \;=\; \bigcup_{q(\tilde\eta)\ge q(\eta)} \CS^H_{\tilde\eta}\,.
\end{equation}
The \emph{null cone} $N_G(V)=\{v\in V\setminus\{0\}\colon 0\in\ol{G\cdot v}\}$ equals the union $\bigcup_{\eta\in\CB_H(A,V)}\CS^H_\eta$. But note that it can occur that $\CS^H_\eta=\varnothing$ for some $\eta\in\CB_H(A,V)$.
\end{definition}

The image of the set $\CB_H(A,V)$ under the map \eqref{eqn:M_to_p} lies in $\Fa/W(\Fg,\Fa)=\Fp/\Ad(K)$. On the other hand, the set $\CB_{BL}$ of stratum labels considered in the beginning is used only through their $\Ad(K)$-orbit. As $\Fp=\bigcup_{g\in K}\Ad_g(\Fa)$, the set $\CB_{BL}$ gives rise to a subset of $\Fp/\Ad(K)$, which we denote $\CB_{BL}(A,V)$.

\begin{theorem} \label{ThmStratum_H-BL}
For every $\beta\in \Fp/\Ad(K)$ with $\CS_\beta\ne\varnothing$ there is a unique $\eta\in \CB_H(A,V)$ with $\Lie(\eta)=\beta$ and $\CS_\beta\;=\; \CS^H_\eta$.
\end{theorem}

In order to prove the theorem we recall from \cite[Definition~7.8]{Boehm-Lafuentec} that the stratum $\CS_\beta$ is described as follows.

\begin{definition}
For $\beta\in\Fp$ set
\begin{equation}\label{eqn:Beta+}
\pi(\beta)^+ \;:=\; \pi(\beta)-\|\beta\|_\Fp^2\cdot\id_V \;\in\; \End_\BR(V)\,.
\end{equation}
\begin{enumerate}
\item
Since $\pi(\beta)^+\in\Sym(V,\sca_{_V})$ by Proposition~\ref{PropExistSkalProd}, it is diagonalizable with real eigenvalues, and $V$ decomposes into orthogonal eigenspaces 
\[
V\;=\;\bigoplus_{r\in\BR} V_{\pi(\beta)^+}^r
\]
$V_{\pi(\beta)^+}^r:=\{v\in V\colon \pi(\beta)^+ (v) = r\cdot v\}=\{v\in V\colon \pi(\beta)(v) = (r+\|\beta\|_\Fp^2)\cdot v\}$. We set
\[
V_{\pi(\beta)^+}^{\ge0} \;:=\;\bigoplus_{r\ge 0} V_{\pi(\beta)^+}^r
\]
and let $p_\beta\colon V_{\pi(\beta)^+}^{\ge 0} \to V_{\pi(\beta)^+}^0$ be the orthogonal projection.
\item Assume that $\beta=\Lie(\eta)$ for $\eta=\tfrac{1}{n}\lambda\in M(G)$ wiwth $\lambda\in X_*(A)$. The centralizer $G_\beta:=Z_G(\lambda)=\{ g\in G\colon \Ad_g(\beta)=\beta\}$ of the cocharacter $\lambda$ is an algebraic subgroup of $G$ and reductive by \cite[\S\,13.17 Corollary~2]{Borel91}. Its Lie algebra $\Fg_\beta$ equals $\{X\in\Fg\colon [\beta,X]=0\}$. We consider $K_\beta:=K\cap G_\beta$ and its Lie algebra $\Fk_\beta$. We define 
\[
\Fh_\beta \;:=\; \Fk_\beta \oplus \{X\in\Fp\colon \langle \beta,X\rangle_\Fp=0\}
\]
and let $H_\beta$ be the algebraic subgroup of $G_\beta$ with Lie algebra $\Fh_\beta$.
\item
The subspace of \emph{$H_\beta$-semi-stable vectors} is defined as
\[
U_{\pi(\beta)^+}^0 \;:=\; \{v\in V_{\pi(\beta)^+}^0\colon 0\notin \ol{H_\beta\cdot v}\,\}\,.
\]
We define $U_{\pi(\beta)^+}^{\ge 0}:=p_\beta^{-1}(U_{\pi(\beta)^+}^0)$ and the \emph{stratum} associated with the orbit $\Ad_K(\beta)$
\begin{equation}\label{eqn:SBeta}
\CS_\beta\;:=\; G\cdot U_{\pi(\beta)^+}^{\ge 0}.
\end{equation}
\end{enumerate}
\end{definition}

Since $\Fp=\bigcup_{h\in K} \Ad_h(\Fa)$ we can conjugate under $K$ and assume that $\beta\in\Fa$. In this case, we restrict the representation $\rho$ to $A$, and hence assume that $G=A$ is a split torus over $\BR$. Then $\langle \phi(\lambda'),\eta\rangle_{\!_A}=\langle\Lie\lambda',\beta\rangle_\Fa=0$ for every $\lambda'\in X_*(H_\beta)$.

\begin{lemma}\label{LemmaCompareStrata}
Assume that $G=A$ is a split torus over $\BR$. Let $\eta=\tfrac{1}{n}\lambda\in M(A)=X_*(A)_\BQ$ with $\lambda\in X_*(A)$ and $\beta=\Lie(\eta)$. Then $q(\eta)=\|\beta\|_\Fa^2$ and
\[
V_{\pi(\beta)^+}^{\ge 0} \;=\; \{v\in V\colon m(v,\tfrac{\eta}{q(\eta)} ) \ge 1\} \qquad\text{and}\qquad V_{\pi(\beta)^+}^r \setminus\{0\} \;\subset\; \{v\in V\colon m(v,\tfrac{\eta}{q(\eta)} ) = 1 + \tfrac{r}{q(\eta)}\}
\]
and
\[
U_{\pi(\beta)^+}^{\ge 0} \;=\; \bigl\{\, v \in V \colon \eta/q(\eta) \in \Lambda_A(v)\,\bigr\}\,.
\]
\end{lemma}

\begin{proof}
We start with the following observation. Let $v=\sum_r v_r\in V$ with $v_r\in V_{\pi(\beta)^+}^r$. Let $v=\sum_\chi v_\chi$ and $v_r=\sum_\chi v_{r,\chi}$ be the decompositions according to weight spaces as in \eqref{eqn:WeigthSpaces}. Then
\begin{equation}\label{eqn:v_rchi}
\langle\chi,\eta\rangle_{\!_A}\cdot v_{r,\chi}=\pi(\beta)\cdot v_{r,\chi}=(r+\|\beta\|_\Fa^2)\cdot v_{r,\chi}\,,
\end{equation}
and hence $\langle\chi,\eta\rangle_{\!_A}=\|\beta\|_\Fa^2$ implies $v_{r,\chi}=0$ for $r\ne 0$ and $v_\chi=v_{0,\chi}$. Therefore, $v_{0,\chi}\ne 0$ if and only if $v_\chi \ne 0$ and $\langle\chi,\eta\rangle_{\!_A}=\|\beta\|_\Fa^2$. Thus, the states satisfy
\begin{equation}\label{eqn:StatesOf_v_vs_v_0}
R(v_0,A) \;=\; R(v,A) \cap \{\chi\in X^*(A)\colon \langle\chi,\eta\rangle_{\!_A} = \|\beta\|_\Fa^2\}
\end{equation}

Now let $v\ne 0$. Then $\lambda(t)\cdot v_r = \exp(n\pi(\beta)\cdot\log t)\cdot v_r= t^{n(r+q(\eta))}\cdot v_r$ for $t\in\BR_{>0}$, and $v_r\ne0$ implies $m(v_r,\lambda)=n(r+q(\eta))$. In particular, $m(v,\tfrac{\eta}{q(\eta)})=\tfrac{1}{n\,q(\eta)}\cdot m(v,\lambda)\ge 1$ if and only if $v_r=0$ for all $r$ with $n(r+q(\eta))<n\,q(\eta)$, that is for all $r<0$. This proves the equality for $V_{\pi(\beta)^+}^{\ge 0}$.

To prove the equality for $U_{\pi(\beta)^+}^{\ge 0}$ let $v=\sum_{r\ge0} v_r\in U_{\pi(\beta)^+}^{\ge 0}$ with $v_r\in V_{\pi(\beta)^+}^r$ and $v_0\in U_{\pi(\beta)^+}^0$. Then $v_0\ne0$, and $m(v_0,\tfrac{\eta}{q(\eta)})=1$ by the above. Let $\Lambda_A(v_0)=\{\tilde\eta\}$ and write $\tilde\beta:=\Lie(\tilde\eta)=\beta'+s\cdot\beta$ for $\beta'\in \Fh_\beta=(\BR\beta)^\perp$ and $s\in\BR$. In fact, $s=\langle\tilde\beta,\beta\rangle_\Fa/\|\beta\|_\Fa^2\in\BQ$, because $\beta,\tilde\beta\in \Lie(X_*(A)_\BQ)$, and hence $\beta'\in \Lie(X_*(A)_\BQ)$.

Let $\eta'=\tfrac{1}{n'}\lambda'$ with $\lambda'\in X_*(A)$ and $\Lie(\eta')=\beta'$. Then $\eta'=\tilde\eta-s\eta$. Moreover, $m(v_0,\tilde\eta)\ge 1=m(v_0,\tfrac{\eta}{q(\eta)})$ implies
\[
\frac{1}{q(\eta)} \;=\; q\bigl(\tfrac{\eta}{q(\eta)}\bigr) \;\ge\; q(\tilde\eta) \;=\; q(\eta')+s^2q(\eta) \;\ge\; s^2q(\eta),
\]
and hence $s\,q(\eta)\le 1$. If $\beta'\ne0$, then $q(\eta')>0$ and we even have a strict inequality $s\,q(\eta)<1$. In this case, $m(v_0,\eta')= m(v_0,\tilde\eta)-s\,q(\eta)>0$ by \eqref{eqn:mOfSum}. By definition of $m(v_0,\lambda')=n'\cdot m(v_0,\eta')>0$ this means that $\lim_{t\to0}\lambda'(t)\cdot v_0=0$. Since $\lambda'(t)\in H_\beta$, we obtain $0\in \ol{H_\beta\cdot v_0}$ in contradiction to $v\in U_{\pi(\beta)^+}^{\ge 0}$. We conclude that $\beta'$ must be $0$ and $\tilde\eta=s\eta$. To show that $\tilde\eta=\eta/q(\eta)$ we observe that
\[
1 \;\le\; m(v_0,\tilde\eta) \;=\; m(v_0,s\eta) \;=\; s\,q(\eta)\cdot m(v_0,\tfrac{\eta}{q(\eta)}) \;=\; s\,q(\eta),
\]
and hence $s\,q(\eta)=1$ as desired.

To show that $\Lambda_A(v)=\{\eta/q(\eta)\}$ we let $\tilde\eta\in X_*(A)_\BQ$ with $\Lambda_A(v)=\{\tilde\eta/q(\tilde\eta)\}$. Since $t^{-q(\eta)}v\in\CS^H_{\tilde\eta}$ for all $t\in\BR\mal$, and $v_0=\lim_{t\to 0}\lambda(t)\cdot t^{-q(\eta)}v$, we see that $v_0\in\ol{\CS^H_{\tilde\eta}}$. Therefore, $q(\eta)\ge q(\tilde\eta)$ by \eqref{eqn:HesslinkStrat}. On the other hand, $m(v,\tfrac{\eta}{q(\eta)})\ge 1$ implies that $\tfrac{1}{q(\tilde\eta)}=q(\tfrac{\tilde\eta}{q(\tilde\eta)})=q^*_A(v)\le q(\tfrac{\eta}{q(\eta)})=\tfrac{1}{q(\eta)}$, and hence $q^*_A(v)=q^*_A(v_0)$. Then $\Lambda_A(v)=\Lambda_A(v_0)=\{\eta/q(\eta)\}$ by \cite[Lemma~2.8a)]{Hesselink79} as desired.

Conversely, let $v=\sum_r v_r\in V$ with $v_r\in V_{\pi(\beta)^+}^r$ and $\Lambda_A(v)=\{\eta/q(\eta)\}$. In particular, $m(v,\tfrac{\eta}{q(\eta)})\ge 1$ and $v\in V_{\pi(\beta)^+}^{\ge 0}$. If $v_0=0$, then $m:=m(v,\tfrac{\eta}{q(\eta)})>1$, and $\tfrac{\eta}{m\,q(\eta)}$ has $m(v,\tfrac{\eta}{m\,q(\eta)})=1$ and $q(\tfrac{\eta}{m\,q(\eta)})= q(\tfrac{\eta}{q(\eta)})/m^2 < q(\tfrac{\eta}{q(\eta)})$ in contradiction to $\Lambda_A(v)=\{\eta/q(\eta)\}$. Thus, $v_0\ne 0$, and we claim that $\Lambda_A(v_0)=\Lambda_A(v)=\{\eta/q(\eta)\}$. Indeed, $\eta$ is the element in the convex hull of $\phi^{-1}R(v,A)$ with minimal norm $q(\eta)$. Therefore, $\langle\chi,\eta\rangle_{\!_A} \ge q(\eta)$ for all $\chi\in R(v,A)$, because otherwise $q(c\,\phi^{-1}(\chi)+(1-c)\eta)<q(\eta)$ for $0<c<2(q(\eta)-\langle\chi,\eta\rangle_{\!_A})/q(\phi^{-1}(\chi)-\eta)$. If now $\eta=\sum_{\chi\in R(v,A)}c_\chi\cdot \phi^{-1}(\chi)$ with $\sum_\chi c_\chi=1$, then $q(\eta) = \sum_\chi c_\chi\langle\chi,\eta\rangle_{\!_A}\ge\sum_\chi c_\chi q(\eta)$ implies that $c_\chi=0$ whenever $\langle\chi,\eta\rangle_{\!_A} >q(\eta)$. By \eqref{eqn:StatesOf_v_vs_v_0} we conclude that $\eta$ lies in the convex hull of $R(v_0,A)$, and hence $\Lambda_A(v_0)=\{\eta/q(\eta)\}$ as claimed.

Assume that $0\in\ol{H_\beta\cdot v_0}$. Then there is a cocharacter $\lambda'\in X_*(H_\beta)$ with $m(v_0,\lambda')\ge 1$ by the Hilbert-Mumford criterion; see for example~\cite[1.2]{Hesselink79}. Let $\eta'$ be a rational multiple of $\lambda'$ with $m(v_0,\eta')=q(\eta)$. Let $0<c<2q(\eta)/(q(\eta')+q(\eta))$ and set $\tilde\eta:=c\eta'+(1-c)\eta$. Then we have $m(v_0,\tilde\eta)=q(\eta)$ by \eqref{eqn:mOfSum}. But $q(\tilde\eta)=c^2q(\eta')+(1-c)^2q(\eta)<q(\eta)$. Therefore, $q\bigl(\tilde\eta/q(\eta)\bigr)<1/q(\eta)$ and $m(v_0,\tilde\eta/q(\eta))=1$ in contradiction to $\Lambda_A(v_0)=\{\eta/q(\eta)\}$. This shows that $0\notin\ol{H_\beta\cdot v_0}$ and $v\in U_{\pi(\beta)^+}^{\ge 0}$ proving the lemma.
\end{proof}

\begin{proof}[Proof of Theorem~\ref{ThmStratum_H-BL}]
For $\eta\in X_*(A)_\BQ$ representing a point of $\CB_H(A,V)$ let $\beta=\Lie(\eta)\in\Fp$. We show that $\CS_\beta=\CS^H_\eta$. Let $v\in \CS^H_\eta$. Then there is an element $g\in G$ with $\eta/q(\eta) \in \Lambda_G(g\cdot v)$. In particular, $q^*_A(g\cdot v)= q^*_G(g\cdot v)$, and this implies that $A$ is optimal for $v$. By Lemma~\ref{LemmaCompareStrata} we have $g\cdot v\in \bigl\{\, \tilde v \in V \colon \eta/q(\eta)\in \Lambda_A(\tilde v)\,\bigr\} = U_{\pi(\beta)^+}^{\ge 0} $, and hence $v\in \CS_\beta$. This implies further
\[
\CS_\beta \;=\; G\cdot U_{\pi(\beta)^+}^{\ge 0} \;=\; G\cdot \bigl\{\, \tilde v \in V \colon \eta/q(\eta) \in \Lambda_G(\tilde v)\,\bigr\} \;=\;\CS^H_\eta\,.
\]
Let now $\beta\in\Fp/\Ad(K)$ be such that $\CS_\beta\ne\varnothing$, and let $v\in\CS_\beta$. Then there is a $g\in G$ such that $A$ is optimal for $g\cdot v$. Let $\eta\in X_*(A)_\BQ$ with $\Lambda_A(g\cdot v)=\{\eta/q(\eta)\}$. Then $\eta\in\CB_H(A,V)$ and $v\in\CS^H_\eta=\CS_{\Lie(\eta)}$ from the above. Therefore $\beta=\Lie(\eta)$ and the theorem is proven.
\end{proof}

\subsection{Changing the group}\label{SubsectChangingGp}

Let $f\colon G'\to G$ be an (algebraic) morphism of real reductive groups. Let $K\subset G$ and $K'\subset G'$ be maximal compact subgroups with $f(K')\subset K$. The canonical Cartan decompositions $\Fg=\Fk\oplus\Fp$ and $\Fg'=\Fk'\oplus\Fp'$ from \S\,\ref{PointGKp} satisfy $Df(\Fp')\subset \Fp$. Let $A'\subset G'$ and $A\subset G$ be maximal split tori as in \S\,\ref{PointAbelianSub} with $\Fa'\subset\Fp'$ and $\Fa\subset\Fp$, such that $f(A')\subset A$. This induces maps $f_*\colon M(G')\to M(G)$ and  $f_*\colon M(A')\to M(A)$ given by $\lambda'\mapsto f\circ\lambda'$, and $f\dual\colon X^*(A)\to X^*(A'),\ \chi\mapsto \chi\circ f$. If $\rho\colon G\to\GL(V)$ is a representation and $\rho':=\rho\circ f\colon G'\to\GL(V)$ then $V_\chi\subset V_{f\dual\chi}$ for every $\chi\in X^*(A)$ and $V_{\chi'}=\bigoplus\limits_{\chi\colon f\dual\chi=\chi'}V_\chi$. In particular, $X^*(A',V)$ is the image of $X^*(A,V)$ under $f\dual$, and $R(x,A')$ is the image of $R(x,A)$ under $f\dual$ for every $x\in V$.

\subsubsection*{Case $f_*\colon M(A')\into M(A)$ is injective.}

Let $M(A')^\perp$ be the orthogonal complement of $M(A')$ in $M(A)$ with respect to the scalar product $(\,.\,,.\,)$ on $M(A)$ associated with the quadratic form $q$. We restrict both to the subspace $M(A')$, i.e.~$q'\colon M(A')\to\BQ$ is given by $q'(\eta'):=q(f_*\eta')$. The isomorphism $\phi\colon M(A)\isoto X^*(A)_\BQ$ is determined by the condition that $(m_1,m_2)=\langle\phi(m_1),m_2\rangle_A$ for all $m_i\in M(A)$. We note that $\langle f\dual\chi,\lambda'\rangle_{A'} = \chi\circ f\circ\lambda' = \langle \chi,f_*\lambda'\rangle_A$ for $\lambda'\in X_*(A')$ and $\chi\in X^*(A)$. Therefore, $\phi':=f\dual\circ \phi\circ f_*\colon M(A')\isoto X^*(A')_\BQ$ is the isomorphism with $\langle \phi'(m'_1),m'_2\rangle_{A'}=\langle \phi(f_* m'_1),f_* m'_2\rangle_A=(f_* m'_1,f_* m'_2)$ for all $m_i'\in M(A')$. In the diagram
\[
\xymatrix{
0 \ar[r] & M(A') \ar[d]^{\phi'} \ar[r]^{f_*} & M(A) \ar[d]^\phi \ar[r]^{pr} & \coker(f_*) \ar[r] & 0 \\
0 & \ar[l] X^*(A')_\BQ & \ar[l]^{f\dual} X^*(A)_\BQ & \ar[l] \ker(f\dual) & \ar[l] 0
}
\]
the orthogonal complement $M(A')^\perp$ is equal to the preimage of $\ker(f\dual)$ under $\phi$, because $0=(m,f_* m')=\langle \phi(m),f_* m'\rangle_A=\langle f\dual\phi(m),m'\rangle_{A'}$ for every $m'\in M(A')$ if and only if $\phi(m)\in \ker(f\dual)$. It follows that the orthogonal projection of $M(A)$ onto $M(A')$ equals the composition $(\phi')^{-1}f\dual\phi$.

If $\rho\colon G\to\GL(V)$ is a representation and $x\in V$, it follows that (the convex hull of) $\phi'{}^{-1}R(x,A')$ is the orthogonal projection of (the convex hull of) $\phi^{-1}R(x,A)$. However, the minimum $\eta(x,A)$ is in general \emph{not} mapped to the minimum $\eta(x,A')$.

\subsubsection*{Case $f_*\colon M(A')\onto M(A)$ is surjective.}

We need a quadratic form $q'$ and scalar product $(\,.\,,.\,)'$ on $M(A')$ such that the orthogonal complement $\ker(f_*)^\perp$ of $\ker(f_*)\subset M(A')$ maps under $f_*$ isometrically onto $M(A)$. If $q'\colon M(A')\to\BQ$ is given, we set $q(\eta):=\min\{q'(\eta')\colon f_*\eta'=\eta\}$ for $\eta\in M(A)$. On the other hand, if $q$ is given, we achieve this by choosing a splitting $s\colon M(A)\into M(A')$ of $f_*$ and a quadratic form $q''$ on $\ker(f_*)$, and setting $q'(\eta)=q(f_*\eta)+q''(\eta-s\,f_*\eta)$. For $\chi\in X^*(A)_\BQ$ and $\eta'\in \ker(f_*)^\perp\subset M(A')$ we then obtain $(f_*\phi'{}^{-1}f\dual\chi,f_*\eta')'=(\phi'{}^{-1}f\dual\chi,\eta')=\langle f\dual\chi,\eta'\rangle_{A'}=\langle \chi,f_*\eta'\rangle_A$, and hence $\phi=(f_*\phi'{}^{-1}f\dual)^{-1}$.
\[
\xymatrix{
0 \ar[r] & \ker(f_*) \ar[r] & M(A') \ar[d]^{\phi'} \ar[r]^{f_*} & M(A) \ar[d]^\phi \ar[r] & 0 \\
0 & \ar[l] \coker(f\dual) & \ar[l] X^*(A')_\BQ & \ar[l]^{f\dual} X^*(A)_\BQ & \ar[l] 0
}
\]
If $\rho\colon G\to\GL(V)$ is a representation and $x\in V$, it follows that $f\dual\colon R(x,A)\isoto R(x,A')$ and $\phi^{-1}R(x,A)=f_*\phi'{}^{-1}R(x,A')$. But again the minimum $\eta(x,A')$ is in general \emph{not} mapped to the minimum $\eta(x,A)$.

\subsubsection*{Hesselink's results.}

Hesselink makes the assumption that $f$ is \emph{special}, by which he means that $f^\der\colon (G')^\der\to G^\der$ is surjective for the derived groups, and $f\dual$ maps the roots of $G$ bijectively onto the roots of $G'$. Moreover, he assumes that $V$ is \emph{of adjoint type}; see \cite[(7.1)~Definition]{Hesselink78}. Under these condition he proves in \cite[(10.4)~Proposition]{Hesselink78}
\begin{itemize}
\item $f_*$ induces a bijection $\Lambda_{G'}(x) \isoto \Lambda_G(x)$ for every $x\in V$.
\item If $x$ is unstable for the action of $G$, then $f^{-1}Q_{G,x}=Q_{G',x}$ for the parabolic subgroups $Q_{G,x}\subset G$ and $Q_{G',x}\subset G'$ from \eqref{eqn:ParabolicQ_beta}.
\end{itemize}

\section{Example: Jordan Normal Form}

Let $G=\GL_n(\BR)$ and $K=\SO_n$. Then $\Fg\Fl_n=\Fs\Fo_n\oplus\Fp$ is the canonical Cartan decomposition with $G=K\cdot\exp\Fp$, where $\Fp=\Sym_n=\{X\in\Fg\Fl_n\colon X=X^\top\}$ and $\Fs\Fo_n=\Lie\SO_n=\{X\in\Fg\Fl_n\colon X+X^\top=0\}$. 

Let $\Fa\subset\Fp$ be the sub-algebra of diagonal matrices. Then $A\subset G$ is the split torus of diagonal matrices. We have $M:=\{g\in G\colon \Ad(g)|_\Fa=\id_\Fa\}=A$ and ${}^{0\!}M:=M\cap K=(1)$, because $\Lie{}^{0\!}M=\Lie M\cap\Fk=\Fa\cap\Fs\Fo_n=(0)$. So $\Fs=(0)$ and $\Ft=\Fa$, whence $T=A$. We fix the scalar product $\langle X,Y\rangle_{\Fg\Fl_n}:=\tr(X\cdot Y^\top)=\sum_{i,j} X_{ij}Y_{ij}$ on $\Fg\Fl_n$, $\Fp$ and $\Fa$. This scalar product on $\Fp$ is $\Ad(K)$-invariant. On $\Fa$ and $\BQ^n=M(A)\subset\Fa=\BR^n$ it induces the standard scalar product. The induced isomorphism $\phi$ equals $M(A)=\BQ^n\xrightarrow{\ \id\;}\BQ^n=X^*(A)_\BQ$. 

\subsection*{The right regular representation.} $G=\GL_n(\BR)$ acts on the right regular representation $V=\BR^n$ by $\rho(g)\colon v\mapsto g\cdot v$ for $g\in G$, and hence $\pi=\id_{\Fg\Fl_n}$. Here $V$ equals the null cone. The weight space decomposition $V=\bigoplus_{\chi\in X^*(A)} V_\chi$ is given by $V_\chi=\BR\cdot e_i$ if $\chi=e_i\colon\diag(t_1,\ldots,t_n)\mapsto t_i$ and $V_\chi=(0)$ otherwise. So $X^*(A,V)=\{e_1,\ldots,e_n\}$ and
\[
\CB_H(A,V)\;=\;\bigl\{\,\tfrac{1}{\#J}\sum_{i\in J}e_i\ \colon\ \varnothing\ne J\subset\{1,\ldots,n\}\,\bigr\}\;=\;\bigl\{\,e_i, \ \tfrac{1}{2}(e_i+e_j),\ \ldots\,\bigr\}\,.
\]
Not every $\eta\in\CB_H(A,V)$ occurs as the label of a non-empty stratum $\CS_{\beta(\eta)}\ne\varnothing$, but only the $\eta=e_i$ do. These all form the orbit of $e_1$ under $W(\Fg,\Fa)=\CS_n$. Indeed, every $0\ne v\in V$ lies in the $G$-orbit of $e_1$, $A$ is optimal for $e_1$, and $R(e_1,A)=\{e_1\}$, whence $\eta(e_1,A)=e_1$. Its image under the map~\eqref{eqn:M_to_p} considered in Theorem~\ref{ThmStratum_H-BL} is the matrix $\beta:=\beta(\eta)=e_1\cdot e_1^\top\in\Fa\subset\Fp$. In particular, up to the action of $K$, there is only one non-empty stratum label in $\CB_{BL}(A,V)$, namely $K\cdot\beta=K\cdot e_1 e_1^\top$. We compute $\|\beta\|^2=\tr(\beta^2)=1$. We then have
\[
\pi(\beta)_+:= \pi( \beta)- \Vert \beta \Vert_\Fp^2\cdot \id_V =\beta-\Id_n = -e_2 e_2^\top - \ldots - e_n e_n^\top\;\in\;\End_\BR(V)
\]
which is negative semi-definite on all of $V$. We also have $V_{\beta_+}^{\ge0}=V_{\beta_+}^0=\BR\cdot e_1\subset V$ on which $\pi(\beta)_+$ is zero and $U_{\beta_+}^{\ge0}=U_{\beta_+}^0=\BR\mal\cdot e_1$. Since $\pi=\id_{\Fg\Fl_n}$ we can take $\beta_+=\beta-\Id_n$ which is a negative semi-definite symmetric matrix in $\Fp$.

To compute the moment map we consider the standard scalar product $\langle v,w\rangle_{_V} =v^\top \cdot w$ on $V=\BR^n$ and let $S \in \Fp$. Then
$$
  \langle m(v),S\rangle_\Fp =\tr (m(v) \cdot S)\overset{!}{=}\frac{\langle \pi(S).v,v\rangle_{_V}}{\Vert v\Vert^2} = \frac{\tr (v^\top S v)}{\Vert v\Vert^2} = \frac{\tr \big((v \cdot v^\top) \cdot S\big)}{\Vert v\Vert^2}\,.
$$
Thus
$$
   m(v)=\frac{v \cdot v^\top}{\Vert v\Vert^2}\,.
$$
This implies that the energy $F(v)=\tr\bigl(m(v)\cdot m(v)\bigr)=1$ is constant, which shows again that there is only one stratum.

\subsection*{The dual of the right regular representation.} $G$ acts on $V\dual:=\Hom_\BR(\BR^n,\BR)$ by $\rho(g)\colon\check v\mapsto \check v\cdot g^{-1}$ for $g\in G$, and hence $\pi(X) = -X^\top$. Again $V\dual$ equals the null cone. The weight space decomposition $V\dual=\bigoplus_{\chi\in X^*(A)} V\dual_\chi$ is given by $V\dual_\chi=\BR\cdot\check e_i$ if $\chi=(-e_i)\colon\diag(t_1,\ldots,t_n)\mapsto t_i^{-1}$ and $V\dual_\chi=(0)$ otherwise. So $X^*(A,V)=\{-e_1,\ldots,-e_n\}$ and
\[
\CB_H(A,V)\;=\;\bigl\{\,\tfrac{1}{\#J}\sum_{i\in J}-e_i\ \colon\ \varnothing\ne J\subset\{1,\ldots,n\}\,\bigr\}\;=\;\bigl\{\,-e_i, \ -\tfrac{1}{2}(e_i+e_j),\ \ldots\,\bigr\}\,.
\]
Also here only the $\eta=-e_i$ occur as the label of a non-empty stratum $\CS_{\beta(\eta)}\ne\varnothing$. Then $\beta:=\beta(\eta)=-e_i e_i^\top$ and $\|\beta\|^2=1$. So again
\[
\pi(\beta)_+:= \pi( \beta)- \Vert \beta \Vert_\Fp^2\cdot \id_V =-\beta-\Id_n = - \sum_{j\ne i}e_j e_j^\top\;\in\;\End_\BR(V\dual)
\]
is negative semi-definite on $V\dual$. We again have $(V\dual)_{\beta_+}^{\ge0}=(V\dual)_{\beta_+}^0=\BR\cdot \check e_i\subset V\dual$ on which $\pi(\beta)_+$ is zero and $U_{\beta_+}^{\ge0}=U_{\beta_+}^0=\BR\mal\cdot\check e_1$. However, this time $\pi(X)=-X^\top$, and hence we can take $\beta_+=-\beta+\Id_n$ which is a \emph{positive} semi-definite symmetric matrix in $\Fp$.

\subsection*{The adjoint representation.} It equals $\Fg\Fl_n=V\otimes V\dual$ with $\rho(g)=\Ad_g\colon x\mapsto gxg^{-1}$ for $g\in G$ and $x\in\Fg\Fl_n$, and $\pi=\ad\colon\Fg\to\End_\BR(\Fg\Fl_n)$. Here the null cone consists of the nilpotent matrices. The weight space decomposition $\Fg\Fl_n=\bigoplus_{\chi\in X^*(A)} \Fg\Fl_{n,\chi}$ is given by $\Fg\Fl_{n,\chi}=\BR\cdot e_i\otimes \check e_j$ if $\chi=e_i-e_j\colon\diag(t_1,\ldots,t_n)\mapsto t_i\cdot t_j^{-1}$ and $\Fg\Fl_{n,\chi}=(0)$ otherwise. So $X^*(A,\Fg\Fl_n)=\{e_i-e_j\colon i,j\}$ is the set of roots together with $0$, and $\CB_H(A,\Fg\Fl_n)$ is a bit complicated to describe.

Also here not every $\eta\in\CB_H(A,\Fg\Fl_n)$ satisfies $\CS_{\beta(\eta)}\ne\varnothing$. Namely, the $G$-orbits of nilpotent matrices in $\Fg\Fl_n$ are given by the Jordan normal form and have a representative $x_J:=\sum_{i\in J}e_i\otimes \check e_{i+1}$ where $J\subset\{1,\ldots,n-1\}$. We write $J=\{1,\ldots,n\}\setminus\{n_1,n_1+n_2,\ldots,\sum_{i=1}^s n_i\}$ with $\sum_{i=1}^s n_i=n$ and $n_1\ge n_2\ge\ldots\ge n_s$ giving the sizes of the Jordan-blocks in $x_J$. Then $R(x_J,A)=\{e_i-e_{i+1}\colon i\in J\}$. (This is a subset of the positive simple roots $\{e_i-e_{i+1}\colon \text{all }i\}$ with respect to the Borel subgroup of upper triangular matrices.) Then
\begin{eqnarray}
\nonumber \eta_J\;:=\;\eta(x_J,A) & = & \Bigl(\frac{n_1-1}{2},\frac{n_1-3}{2},\ldots,\frac{1-n_1}{2},\quad \frac{n_2-1}{2},\ldots,\frac{1-n_2}{2},\quad \fdot\ \fdot\ \fdot\ \Bigr) \qquad\text{with}\\[2mm]
\nonumber \beta_J\;:=\;\beta(\eta_J) & = & \diag\Bigl(\frac{n_1-1}{2},\frac{n_1-3}{2},\ldots,\frac{1-n_1}{2},\quad \frac{n_2-1}{2},\ldots,\frac{1-n_2}{2},\quad \fdot\ \fdot\ \fdot\ \Bigr)\;\in\;\Fa\,.
\end{eqnarray}
Note that for $4\le n_1\le n$ this $\eta_J$ does not lie in $\CB_H(A,V)+\CB_H(A,V\dual)$. In any case we will see below that
\begin{equation}\label{eqn:For_q_OfBeta_JordanNF}
q(\eta_J) \;=\;\|\beta_J\|^2 \; = \; \sum_{j=1}^s \frac{(n_j-1)\,n_j\,(n_j+1)}{12}\,. 
\end{equation}
Then $n_j\ge 2$ implies $6\le n_j(n_j+1)$, and so for any $n_j\ge 1$ we have
\begin{equation}\label{eqn:JordanNegDef}
\frac{n_j-1}{2}\;\le\;\frac{(n_j-1)\,n_j\,(n_j+1)}{12}\;\le\;\|\beta_J\|^2.
\end{equation}
The endomorphism $\pi(\beta_J)_+:= \pi( \beta_J)- \Vert \beta_J \Vert^2\cdot \id_{\Fg\Fl_n}\in\End_\BR(\Fg\Fl_n)$ is given by $x\mapsto [\beta_J,x]-\|\beta_J\|^2\cdot x$ for $x\in\Fg\Fl_n$. It has a basis of eigen vectors given by $e_k\otimes\check e_l$ for $1\le k,l\le n$, and equations~\eqref{eqn:For_q_OfBeta_JordanNF} and \eqref{eqn:JordanNegDef} show that again $\pi(\beta_J)_+$ is negative semi-definite on $\Fg\Fl_n$. 

To prove formula \eqref{eqn:For_q_OfBeta_JordanNF} we observe that for $x\in\BR$ and $k\in\BN$ 
\[
\sum_{i=0}^{k-1}(i+x)^2\;=\;k\cdot x^2+\sum_{i=0}^{k-1}2i \cdot x +\sum_{i=0}^{k-1} i^2 \;=\; k\cdot x^2 + k(k-1)\cdot x + \tfrac{(k-1)\,k\,(2k-1)}{6}
\]
So if $n=2k$ is even we have
\begin{eqnarray*}
\Bigl(\frac{n-1}{2}\Bigr)^2+\ldots+\Bigl(\frac{1-n}{2}\Bigr)^2 & = & 2\sum_{i=0}^{k-1}(i+\tfrac{1}{2})^2 \\[2mm]
& = & 2\cdot\Bigl(k\cdot\tfrac{1}{4}+\tfrac{k(k-1)}{2}+ \tfrac{(k-1)\,k\,(2k-1)}{6}\Bigr) \\[2mm]
& = & \frac{(n-1)\,n\,(n+1)}{12}\;,
\end{eqnarray*}
and if $n=2k+1$ is odd we have
\begin{eqnarray*}
\Bigl(\frac{n-1}{2}\Bigr)^2+\ldots+\Bigl(\frac{1-n}{2}\Bigr)^2 & = & 2\sum_{i=0}^{k}i^2 \\[2mm]
& = & \frac{k\,(k+1)(2k+1)}{3} \\[2mm]
& = & \frac{(n-1)\,n\,(n+1)}{12}\;.
\end{eqnarray*}
We compute again the moment map for this action. For $S \in \Fp$ we have
$$
 \langle \pi(S).x,x\rangle_{_V} =\tr \big( [S,x] \cdot x^\top\big)=\tr \big( S \cdot [x,x^\top]\big)
$$
thus
$$
   m(x)=\frac{[x,x^\top]}{\Vert x\Vert^2}\,.
$$
Thus we have $m(x)=0$, if and only if $x$ is a normal matrix.

\subsection*{The action on $\Lambda^2(\BR^n)$}

We consider the action of $\GL_n$ on $\Lambda^2(\BR^n)$ given by
$$
\rho(g)(x \wedge y) =g.x \wedge g.y
$$
for $g\in \GL_n$ and $x\wedge y\in\Lambda^2(\BR^n)$. The standard scalar product on $\BR^n$ induces a canonical scalar product on $\Lambda^2(\BR^n)$.
We will use the isomorphism $\Phi:\Lambda^2(\BR^n) \isoto \Fs\Fo_n=\{A\in\Fg\Fl_n\colon A^\top=-A\}\subset\Fg\Fl_n=\End_\BR(\BR^n)$ given by
$$
\Phi(x\wedge y)\;:=\; xy^\top-yx^\top,\qquad \Phi(x \wedge y)(z)\;:=\;\;x\,y^\top z-y\,x^\top z\;=\; x \langle y,z\rangle_{_V} -y\langle x,z\rangle_{_V}\;\,,
$$
where $z \in \BR^n$.
To make $\Phi$ an isometry we use on  $\Fs\Fo_n$ the scalar product
$\langle A,B\rangle_{\Fs\Fo_n}=-\tfrac{1}{2}\tr (AB)$. Now
$$
 \Phi( g\cdot x \wedge g\cdot y  )= g \cdot \Phi(x\wedge y) \cdot g^\top\,. 
 $$
Thus the representation $\rho$ is isomorphic to the representation of $\GL_n$ on $\Fs\Fo_n$ given by $\rho(g)(A)=g\cdot A\cdot g^\top$ with linearization $\pi(g)=D\rho(g)\colon A\mapsto gA+Ag^\top$.
We compute again the moment map for $A \in \Fs\Fo_n$ and $S\in \Fp$:
$$
  \langle \pi(S).A,A\rangle_{\Fs\Fo_n} =-\tfrac{1}{2}\cdot \tr \big( (SA + AS)\cdot A\big)=-\tr S \cdot A^2\,.
$$
Thus
$$
  m(A)=-\frac{A^2}{\Vert A\Vert^2}=\frac{A^\top A}{\Vert A\Vert^2} \geq 0
$$
with energy
$$
 F(A)=\frac{\tr (A^4)}{(\tr A^2)^2}\,.
$$

\subsection*{Changing the group}

The adjoint representation is of adjoint type, and the morphism $f\colon\SL_n\into\GL_n$ is special; see Subsection~\ref{SubsectChangingGp}. If $A':=A\cap\SL_n\subset\SL_n$ is the split torus of diagonal matrices, then $f_*\colon M(A')\into M(A),\ \lambda\mapsto f\circ\lambda$ and $f_*\colon \Fa':=\Fa\cap\Fs\Fl_n\into\Fa$ induce bijections $\Lambda_H(x) \isoto \Lambda_G(x)$ for every $x\in \Fg\Fl_n$, and map the stratum label $\eta_{\SL_n}(x_J,A')\in\Fa'/W(\Fs\Fl_n,\Fa')$ to $\eta_{\GL_n}(x_J,A)\in\Fa/W(\Fg\Fl_n,\Fa)$ for every nilpotent matrix $x_J$ as above. Note that $W(\Fs\Fl_n,\Fa')=W(\Fg\Fl_n,\Fa)$.

For the representations $V$ and $V\dual$ this is not true, as they are not of adjoint type. We have
\[
M(A')\;=\;\bigl\{\, \eta=(\eta_1,\ldots,\eta_n)\colon\sum_i \eta_i=0\,\bigr\} \qquad\text{and}\qquad M(A')^\perp\;=\;\BQ\cdot (1,\ldots,1)
\]
under the identifications $A=\BG_m^n$ and $M(A)=\BQ^n$. So for $e_i\in V=\BR^n$ we have $\eta(e_i,A)=e_i$ and $\eta(e_i,A')= e_i-\tfrac{1}{n}(1,\ldots,1)$. On the other hand,  for $\check e_i\in V\dual$ we have $\eta(\check e_i,A)=-e_i$ and $\eta(\check e_i,A')= -e_i+\tfrac{1}{n}(1,\ldots,1)$.

\section{Example: The Space of Brackets}

We consider the representation of $\GL_n$ on $V=\Lambda^2(\BR^n)\dual\otimes\BR^n=\Hom_\BR(\Lambda^2(\BR^n),\BR^n)$ given by base change:
$(\rho(g).\mu)(x,y)=g\mu(g^{-1}x,g^{-1}y)$ for $g\in\GL_n$ and $\mu\in V$. Let $A\subset\GL_n$ be the diagonal torus. Then $X^*(A,\Lambda^2(\BR^n)\dual) = \{-e_i-e_j\colon 1\le i<j\le n\}$ and $X^*(A,\BR^n)=\{e_l\colon 1\le l\le n\}$, and hence
\[
X^*(A,V)=X^*(A,\Lambda^2(\BR^n)\dual)+X^*(A,\BR^n)=\{e_l-e_i-e_j\colon 1\le l\le n, 1\le i<j\le n\}\,.
\]
On $V$ we fix the scalar product given by
\[
\langle \mu,\mu'\rangle_{_V} \;:=\;\sum_{i,j} \bigl\langle \mu(e_i,e_j),\mu'(e_i,e_j)\bigr\rangle_{\BR^n}\qquad\text{for}\quad \mu,\mu'\in V\,.
\]

We compute the moment map. We have
for a symmetric $(n\times n)$-matrix $S$ and a bracket $\mu\in V$ with $\Vert \mu\Vert_V^2=1$ 
$$
 \tr(m(\mu) \cdot S)=\langle m(\mu),S\rangle_\Fp =\langle \pi(S).\mu,\mu \rangle_{_V}
$$
with
$$
 \pi(S).\mu(x,y)=S\mu(x,y)-\mu(Sx,y)-\mu(x,Sy)
$$
for all $x,y \in \BR^n$.

Let now $S:= e_1 e_1^\top$, i.e.~$S(e_1)=e_1$ and $S(e_i)=0$ for $i\ge 2$. Then
\begin{eqnarray*}
 \langle \pi(S).\mu,\mu \rangle_{_V} 
   &=&
    \sum_{i,j} \langle (\pi(S).\mu)(e_i,e_j),\mu(e_i,e_j)\rangle_{\BR^n} \\
      &=&
  \sum_{i,j} \langle  S\mu(e_i,e_j)-\mu(Se_i,e_j)-\mu(e_i,Se_j)  ,\mu(e_i,e_j)\rangle_{\BR^n} \\
    &=&
      \sum_{i,j} \langle  S\mu(e_i,e_j) ,\mu(e_i,e_j)\rangle_{\BR^n} 
     -2   \sum_{i,j} \langle \mu(Se_i,e_j)  ,\mu(e_i,e_j)\rangle_{\BR^n} \\
     &=&
      \sum_{i,j} \langle \mu(e_i,e_j),e_1\rangle_{_V}^2  
     -2   \sum_{j} \Vert \mu(e_1,e_j)\Vert^2  
\end{eqnarray*}
This shows that
$$
  \langle m(\mu)x,x\rangle_{_V}=  \sum_{i,j} \langle \mu(e_i,e_j),x\rangle_{_V}^2  
     -  2 \sum_{j} \Vert \mu(x,e_j)\Vert^2  \,.
$$
for all $x \in \BR^n$. 

We consider now the case that $\mu$ is a nilpotent Lie bracket, with 
$\mu \in \CS_{\beta_\Fn}$. 
By general theory we may assume that $m(\mu)=\beta_\Fn$ and $\mu \in U^0_{\beta_\Fn^+}$; have to go to the closure of $G.\mu$ in general; $\mu$ is then a critical point for the energy.
We set
$$
  \beta^+_\Fn:= \beta_\Fn + \Vert \beta_\Fn\Vert^2 \cdot \Id_n \in \Fp
$$
Then we have $\pi(\beta^+_\Fn)=\pi(\beta_\Fn)-\Vert \beta_\Fn\Vert^2 \cdot \Id_V$ as in 
\cite{Boehm-Lafuentec} (where we consider $G$ as a subgroup of ${\rm Gl}(V))$.
Using that $\tr \beta_\Fn=-1$ we have $\langle \beta^+_\Fn,\beta_\Fn\rangle_\Fp=\tr(\beta_\Fn\cdot\beta_\Fn)+\tr(\|\beta_\Fn\|^2\cdot \beta_\Fn)= 0$. Thus
for such a $\mu$ with $\Vert \mu\Vert =1$ we have
$$
 0 =\langle m(\mu),\beta^+_\Fn\rangle_\Fp =\langle \pi(\beta^+_\Fn).\mu,\mu\rangle_{_V}
$$
Since $\mu \in U^0_{\beta_\Fn^+}$ we also have
 $D:=\beta^+_\Fn \in {\rm Der}(\mu)$ \footnote{This means $D([x,y])=[x,Dy]+[Dx,y]$.}. Let now $X\in \BR^n$ be an eigen vector with $DX =d\cdot X$.
 We need to show $d>0$. Using $\ad_\mu(X)(Y):=\mu(X,Y)$ we have
 $$
  d\cdot \ad_\mu(X)=\ad_\mu(DX)=[D,\ad_\mu(X)]=[\beta_\Fn,\ad_\mu(X)]\,.
 $$
As a consequence, using $m(\mu)=\beta_\Fn$
and $\tr( [A,B]C)=\tr(A[B,C])$, we obtain
\begin{eqnarray*}
d \cdot \tr \big( \ad_\mu(X)\cdot \ad_\mu(X)^\top)
 &=&
    \tr  \big(  [\beta_\Fn,\ad_\mu(X)]\cdot \ad_\mu(X)^\top \big) \\
    &=&
    \tr \big( \beta_\Fn \cdot [\ad_\mu(X), \ad_\mu(X)^\top] \big) \\
    &=&
     \langle \pi( [\ad_\mu(X), \ad_\mu(X)^\top]).\mu,\mu\rangle_{_V} \\
        &=&
     \big\langle \big(\pi( \ad_\mu(X)) \pi( \ad_\mu(X)^\top)-  \pi( \ad_\mu(X)^\top) \pi( \ad_\mu(X)) \big).\mu,\mu \big\rangle_{_V}
\end{eqnarray*}
Then since by the Jacobi identity $\pi(\ad_\mu(X)).\mu=0$,
this is a quadratic constraint !!!!,  we arrive at
$$
   \langle (\pi( \ad_\mu(X))(\pi( \ad_\mu(X))^\top.\mu. \mu \big\rangle_{_V}
  = \Vert  \big(\pi( \ad_\mu(X)^\top).\mu\Vert^2 \geq 0\,,
$$
using that $\pi(A)^\top=-\pi(A)=\pi(A^\top)$ and $\pi(S)^\top=\pi(S^\top)$ for
skew-symmetric $A$ and symmetric $S$. We deduce $d \geq 0$ provided
that $\ad_\mu(X)\neq 0$.

Suppose now $\ad_\mu(X)=0$. Then on the one hand side
$$
\langle \beta_\Fn X,X\rangle_{_V} = (d- \Vert \beta_\Fn\Vert^2)\cdot \Vert X\Vert^2
$$
but on the other hand side, using $\beta_\Fn=m(\mu)$ and  $\ad_\mu(X)=0$, we see that
that
$$
 \langle \beta_\Fn X,X\rangle_{_V} = \sum_{i,j}\langle \mu(e_i,e_j),X \rangle_{_V}^2 \geq 0\,.
$$
Thus $d>0$.

So we are left we the case $\ad_\mu(X)\neq 0$ and $d=0$. But then, by the above,
$\ad_\mu(X)^\top\in {\rm Der}(\mu)$ as well.
We denote now by $W_1$ the orthogonal complement of $[\Fn,\Fn]$ in $\Fn$,
be $W_2$ the orthogonal complement of $[\Fn,[\Fn,\Fn]]$ in $[\Fn,\Fn]$ and so on.
Then
$$
  \Fn = W_1 \oplus W_2 \oplus \cdots \oplus W_r
$$
Note that any derivation $D$ of $\Fn$ respects $V_i:=W_i \oplus \cdots \oplus W_r$ for
all $i=1,...,r$. This implies, for instance, that
$$
  0=\langle D V_2,W_1 \rangle_{_V} =\langle V_2, D^\top W_1 \rangle_{_V}\,.
$$
Thus $D^\top$ preserves $W_1$,...., $W_r$. For $D=\ad_\mu(X)\neq 0$ this is impossible.

%
%


{\small
\begin{thebibliography}{BGR84}

\bibitem[AH-LH19]{AH-LH} J.~Alper, D.~Halpern-Leistner, J.~Heinloth: \emph{Cartan–Iwahori–Matsumoto Decompositions for Reductive Groups}, preprint 2019 available as 

\bibitem[BL18a]{Boehm-Lafuentea} Christoph~B\"ohm and Ramiro~Lafuente. \emph{Immortal homogeneous Ricci flows}, Invent.\ Math.~{\bfseries 212} (2018), 461--529.

\bibitem[BL18b]{Boehm-Lafuenteb} Christoph~B\"ohm and Ramiro~Lafuente.
  \emph{Homogeneous Einstein metrics on Euclidean spaces are Einstein solvmanifolds},
   arXiv:1811.12594.
	
\bibitem[BL18c]{Boehm-Lafuentec} Christoph~B\"ohm and Ramiro~Lafuente.	
   \emph{Real geometric invariant theory},
   arXiv:1701.00643 (2017).

\bibitem[BH-C62]{B+HC62} A.~Borel, Harish-Chandra: \emph{Arithmetic subgroups of algebraic groups}, Annals of Math.\ (2) {\bfseries 75} (1962), 485--535.

\bibitem[Bor91]{Borel91}A.~Borel: \emph{Linear algebraic groups}, Second Enlarged Edition, Springer-Verlag, Berlin, 1991. 

\bibitem[Del82]{DeligneHodgeCycles} P.~Deligne: \emph{Hodge cycles on abelian varieties}, in ``Hodge Cycles, Motives, and Shimura Varieties'', pp.~1--100, LNM~{\bfseries 900}, Springer-Verlag, New York 1982; also available at \href{https://jmilne.org/math/Documents/Deligne82.pdf}{http:/\hspace{-0.2em}/www.jmilne.org/math}.

\bibitem[HP18]{HartlPal1} Urs~Hartl, Ambrus~P\'al: \emph{Crystalline Chebotar\"ev Density Theorems}, Preprint November 2018 on 
{arXiv:1811.07084}.

\bibitem[HSS08]{HSS08} Peter~Heinzner, Gerald~Schwarz, Henrik~St\"otzel: \emph{Stratifications with respect to actions of real reductive groups}, Compos.\ Math.\ {\bfseries 144} (2008), no. 1, 163--185.

\bibitem[Hes78]{Hesselink78} W.~Hesselink: \emph{Uniform instability in reductive groups}. J.~Reine Angew.\ Math.\ {\bfseries 303/304} (1978), 74--96.

\bibitem[Hes79]{Hesselink79} W.~Hesselink: \emph{Desingularizations of Varieties of Nullforms}. Inventiones Math.\ {\bfseries 55} (1979), 141--163.

\bibitem[Kir84]{Kirwan84} Frances Clare Kirwan: \emph{Cohomology of quotients in symplectic and algebraic geometry}, Mathematical Notes, vol. 31, Princeton University Press, Princeton, NJ, 1984.

\bibitem[Lau10]{Lauret10} Jorge~Lauret: \emph{Einstein solvmanifolds are standard}, Ann.\ of Math.~(2) {\bfseries 172} (2010), no. 3, 1859--1877.

\bibitem[MFK94]{GIT} D.~Mumford, J.~Fogarty, F.~Kirwan: \emph{Geometric invariant theory}, Ergebnisse der Mathematik und ihrer Grenzgebiete (3)~{\bfseries 34}, Springer-Verlag, Berlin, 1994.

\bibitem[Nes84]{Ness84} Linda~Ness: \emph{A stratification of the null cone via the moment map}, Amer.~J.\ Math. {\bfseries 106} (1984), no. 6, 1281--1329, with an appendix by David Mumford.

\bibitem[Ser92]{Serre92} J.-P.~Serre, \emph{Lie algebras and Lie groups}, Lecture Notes in Mathematics {\bfseries 1500}, second edition, Springer-Verlag, Berlin (1992).

\bibitem[Wal88]{Wallach} N.~Wallach: \emph{Real reductive groups I}, Pure and Applied Mathematics, 132. Academic Press, Inc., Boston, MA, 1988.

\end{thebibliography}
}
\end{document}